\documentclass[times,sort&compress,3p]{elsarticle}
\journal{Journal of Multivariate Analysis}
\usepackage[labelfont=bf]{caption}

\usepackage{amsmath,amsfonts,amssymb,amsthm,booktabs,xcolor,epsfig,graphicx,hyperref,url}
\usepackage{MnSymbol}
\usepackage{csquotes}

\def \P{\text{\rm P}}

\DeclareMathOperator*{\Var}{Var}

\newcommand{\inp}[3]{\left\langle #1, #2\right\rangle_{#3}}
\newcommand{\rfrac}[2]{^{#1}\!/_{#2}} 
\definecolor{violet}{rgb}{0.7,0,0.6}
\makeatletter
\newcommand{\labitem}[2]{%
	\def\@itemlabel{\textbf{#1}}
	\item
	\def\@currentlabel{#1}\label{#2}}
\makeatother
\theoremstyle{plain}
\newtheorem{Th}{Theorem}

\newtheorem{Lema}{Lemma}
\newtheorem{Cor}{Corollary}

\theoremstyle{definition}

\begin{document}

	\begin{frontmatter}
		\title{A uniform kernel trick for high and infinite-dimensional two-sample problems}
		
		\author[1]{Javier Cárcamo \corref{mycorrespondingauthor}}
		\author[2]{Antonio Cuevas}
		\author[3]{Luis-Alberto Rodríguez}
		
		\address[1]{Euskal Herriko Unibertsitatea-Universidad del País Vasco, Spain}
		\address[2]{Universidad Autónoma de Madrid and Instituto de Ciencias Matemáticas ICMAT (CSIC-UAM-UCM-UC3M), Spain}
		\address[3]{Georg-August-Universität G\"{o}ttingen, Germany}
		
		\cortext[mycorrespondingauthor]{Corresponding author. Email address: \url{javier.carcamo@ehu.eus}}
		
		\begin{abstract}
			We use a suitable version of the so-called \enquote{kernel trick} to devise two-sample tests, especially focussed on high-dimensional and functional data. Our proposal entails a simplification of the practical problem of selecting an appropriate kernel function. Specifically, we apply a uniform variant of the kernel trick which involves the supremum within a class of kernel-based distances. We obtain the asymptotic distribution of the test statistic under the null and alternative hypotheses. The proofs rely on empirical processes theory, combined with the delta method and Hadamard directional differentiability techniques, and functional Karhunen-Loève-type expansions of the underlying processes. This methodology has some advantages over other standard approaches in the literature. We also give some experimental insight into the performance of our proposal compared to other kernel-based approaches (the original proposal by \cite{Borgwardt2006} and some variants based on splitting methods) as well as tests based on energy distances \cite{Rizzo&Szekely2017}.
			
		\end{abstract}
		
		
	\end{frontmatter}
	\section{Introduction: an overview}\label{Section_introduction}
		In this section we provide an extended summary including not only the main ideas of this work but, specially, the general setting, motivation and related literature, as well as the technical tools we use.
			\subsubsection*{The kernel trick and some potential kernel traps}
				We focus on statistical problems where, essentially, the aim is to properly separate data coming from two different populations; this is the case of binary supervised classification and two-sample testing problems. In such situations, the \textit{kernel trick} is a common paradigm. In a few words, the standard multivariate version (i.e., with data in $\mathbb{R}^{d}$) of the kernel trick lies in separating the data in both populations using a symmetric non-negative definite \enquote{kernel function}. The values of the kernel can be seen as the inner product of transformed versions of the original observations in a different (usually higher-dimensional) space. It is expected that the groups can be better distinguished in the new final space; see \cite{Scholkopf&Smola2018}.

				We are particularly interested in those situations in which the available data are high-dimensional or even functional (thus, infinite-dimensional). In such cases, the strategy of mapping the data into a higher-dimensional space does not seem to be so compelling. Still, the kernel trick remains meaningful in a sort of \enquote{second generation} version, whose point is to take the data to a more comfortable and flexible space. In this new space, the statistical metho\-dology might be mathematically more tractable, and more easily implemented and interpreted. To be more precise, a probability distribution $\operatorname{P}$ on the sample space $\mathcal{X}$ is replaced with a function
				\begin{equation}\label{eq:muP}
					\mu_{\operatorname{P}}(x)=\int_{\mathcal{X}} k(x,y)\, \operatorname{dP}(y),\quad x\in\mathcal{X},
				\end{equation}
				in an appropriate space of \enquote{nice functions} defined by means of the kernel $k$. In this way, the distance between two probability measures is computed in terms of the metric in the functional space. As a matter of fact, one of the most appealing proposals in this direction relies on kernel-based distances, expressed in terms of the embedding transformation $\mu_{\operatorname{P}}$ in \eqref{eq:muP}; see \cite{Borgwardt2006}.

				The kernel $k$ involved in this methodology depends, almost unavoidably, on some tuning parameter $\lambda$, typically a scale factor. Therefore, we actually have a \textit{family} of kernels, $k_{\lambda}$, for $\lambda\in\Lambda$, where $\Lambda$ is usually a subset of $\mathbb{R}^{k}$ ($k\geq1$). For instance, the popular family of \textit{Gaussian kernels} with parameter $\lambda\in(0,\infty)$ is defined by
				\begin{equation}\label{eq:gk}
					k_{\lambda}(x,y)=\operatorname{exp}\left(-\lambda\,\|x-y\|^{2}\right),\quad\operatorname{for}\ x,y\in\mathcal{X},
				\end{equation}
				where $\|\cdot\|$ is a norm in $\mathcal{X}$. Unfortunately, there is no general rule to know \textit{a priori} which kernel works best with the available data. In other words, the choice of $\lambda$ is, to some extent, arbitrary but not irrelevant, as it could remarkably affect the final output. For example, very small or very large choices of  $\lambda$ in \eqref{eq:gk} result in null discrepancies, which have no ability to distinguish distributions. The selection of $\lambda$ is hence a delicate problem that has not been satisfactorily solved so far. This is what we call the \textit{kernel trap}: a bad choice of the parameter leading to poor results. Although this problem was not explicitly considered in \cite{Borgwardt2006} and subsequent works on this topic, the authors were aware of this relevant question; in practice, they use a heuristic choice of $\lambda$.


				Further, a parameter-dependent method might be an obstacle for practitioners who are often reluctant to use procedures depending on auxiliary, hard-to-interpret parameters. We thus find here a particular instance of the trade-off between power and applicability: as stated in \cite{Tukey1959}, the \textit{practical power} of a statistical procedure is defined as \textit{“the product of the mathematical power by the probability that the procedure will be used”} (Tukey credits to Churchill Eisenhart for this idea). From this perspective, our proposal can be viewed as an attempt to make kernel-based homogeneity tests more usable by getting rid of the tuning parameter(s). Roughly speaking, the idea that we propose to avoid selecting a specific value of $\lambda$ within the family  $\left\{k_{\lambda}:\lambda\in\Lambda\right\}$  is to take the supremum over the set of parameters $\Lambda$ of the resulting family of kernel-distances. We call this approach the \textit{uniform kernel trick}, as we map the data into many functional spaces at the same time and use, as test statistic, the supremum of the corresponding kernel distances.  We believe that this methodology could be successfully applied as well in supervised classification, though this topic is not considered in this work.
			\subsubsection*{The topic of this paper}
				Two-sample tests, also called homogeneity tests, aim to decide whether or not it can be accepted that two random elements have the same distribution, using the information provided by two independent samples. This problem is omnipresent in practice on account of their applicability to a great variety of situations, ranging from biomedicine to quality control. Since the classical Student's t-tests or rank-based (Mann-Whitney, Wilcoxon, \dots) procedures, the subject has received an almost permanent attention from the statistical community. In this work we focus on two-sample tests valid, under broad assumptions, for general settings in which the data are drawn from two random elements $X$ and $Y$ taking values in a general space $\mathcal{X}$. The set $\mathcal{X}$ is the \textit{\enquote{sample space}} or \textit{\enquote{feature space}} in the Machine Learning language. In the important particular case $\mathcal{X}=\operatorname{L}^2([0,1])$, $X$ and $Y$ are stochastic processes and the two-sample problem lies within the framework of Functional Data Analysis (FDA).

				Many important statistical methods, including goodness of fit and homogeneity tests, are based on an appropriate metric (or discrepancy measure) that allows groups or distributions to be distinguished. Probability distances or semi-distances reveal to the practitioner the dissimilarity between two random quantities. Therefore, the estimation of a suitable distance helps detect significant differences between two populations. Some well-known, classic examples of such metrics are the Kolmogorov distance, that leads to the popular Kolmogorov-Smirnov statistic, and $\operatorname{L}^{2}$-based discrepancy measures, leading to Cramér-von Mises or Anderson-Darling statistics. These methods, based on cumulative distribution functions, are no longer useful with high-dimensional or non-Euclidean data, as in FDA problems. For this reason we follow a different strategy based on more adaptable metrics between general probability measures.

				The \textit{energy distance} (see the review by \cite{Rizzo&Szekely2017}) and the associated \textit{distance covariance}, as well as \textit{kernel distance}, represent a step forward in this direction since they can be calculated with relative ease for high-dimensional distributions. In \cite{Fukumizu2013} the relationships among these metrics in the context of hypothesis testing are discussed. In this paper we consider an extension, as well as an alternative mathematical approach, for the two-sample test in \cite{Borgwardt2006}. These authors show that kernel-based procedures perform better than other more classical approaches when dimension grows, although they are strongly dependent on the choice of the kernel parameter.
				\subsubsection*{Three important auxiliary tools: RKHS, mean embeddings, and kernel distances}
					To present the contributions of this paper, we briefly refer to some important, mutually related, technical notions. As emphasized in \cite{Berlinet&Thomas-Agnan2011}, \emph{Reproducing Kernel Hilbert Spaces} (RKHS in short) provide an excellent environment to construct helpful transformations in several statistical problems. Given a topological space $\mathcal{X}$ (in many applications $\mathcal{X}$ is a subset of a Hilbert space), a \textit{kernel} $k$ is a real non-negative semidefinite symmetric function on $\mathcal{X}\times\mathcal{X}$. The RKHS associated with $k$, denoted in the following by $\mathcal{H}_{k}$, is the Hilbert space generated by finite linear combinations of type $\sum_{j}\alpha_{j}\,k\left(x_{j},\cdot\right)$; see Section \ref{Section_preliminaries} for additional details.
				
					Let $\mathcal{M}_{\operatorname{p}}(\mathcal{X})$ be the set of (Borel) probability measures on $\mathcal{X}$. Under mild assumptions on $k$, the functions in $\mathcal{H}_{k}$ are measurable and $\operatorname{P}$-integrable, for each $\operatorname{P}\in\mathcal{M}_{\operatorname{p}}(\mathcal{X})$. Moreover, it can be checked that the function $\mu_{\operatorname{P}}$ in \eqref{eq:muP} belongs to $\mathcal{H}_{k}$. The transformation $\operatorname{P}\mapsto\mu_{\operatorname{P}}$ from $\mathcal{M}_{\operatorname{p}}(\mathcal{X})$ to $\mathcal{H}_{k}$ is called the \textit{(kernel) mean embedding}; see \cite{Fukumizu2013} and \cite[Chapter 4]{Berlinet&Thomas-Agnan2011}.  The mean embedding of $\operatorname{P}$ can be viewed as a smoothed version of the distribution of $\operatorname{P}$ through the kernel $k$ within the RKHS. This is evident when $\operatorname{P}$ is absolutely continuous with density $f$ and $k(x,y)=K(x-y)$, for some real function $K$. In this situation, $\mu_{\operatorname{P}}$ is the convolution of $f$ and $K$. On the other hand, mean embeddings appear, under the name of \textit{potential functions}, in some other mathematical fields (such as functional analysis); see \cite[p. 15]{El-Fallah2014}.
					
					The \textit{kernel distance} between $\operatorname{P}$ and $\operatorname{Q}$ in $\mathcal{M}_{\operatorname{p}}(\mathcal{X})$ is
					\begin{equation}\label{definition.kernel.distance}
						d_{k}(\operatorname{P},\operatorname{Q})=\left\|\mu_{\operatorname{P}}-\mu_{\operatorname{Q}}\right\|_{\mathcal{H}_{k}}=\left(\int_{\mathcal{X}^{2}}k\,\operatorname{d}(\operatorname{P}-\operatorname{Q})\otimes(\operatorname{P}-\operatorname{Q})\right)^{1/2},	
					\end{equation}
					where $\|\cdot\|_{\mathcal{H}_{k}}$ stands for the norm in $\mathcal{H}_{k}$ and $(\operatorname{P}-\operatorname{Q})\otimes(\operatorname{P}-\operatorname{Q})$ denotes the product (signed) measure on $\mathcal{X}^{2}=\mathcal{X}\times\mathcal{X}$. Therefore, $d_{k}(\operatorname{P},\operatorname{Q})$ is the RKHS distance between the mean embeddings of the corresponding probability measures. Kernel distances were popu\-la\-rized in machine learning as tools to tackle several relevant statistical problems, such as homogeneity tests \cite{Borgwardt2006}, independence \cite{Fukumizu2007}, test of conditional independence \cite{Fukumizu2008} and density estimation \cite{Fukumizu2011}. The key idea behind this methodology can be seen as a particular case of the fruitful kernel trick paradigm.
				\subsubsection*{Our contributions: the uniform kernel trick}
					We consider a family of kernels $\left\{k_{\lambda}:\lambda\in\Lambda\right\}$,  where $\Lambda$ is certain parametric space. For the Gaussian kernel in \eqref{eq:gk}, $\Lambda=(0,\infty)$, but in general $\lambda$ could be a multidimensional parameter, as in the case of Matérn kernels or inverse quadratic kernels; see \cite [p. 1846]{Sriperumbudur2016}. Each $k_{\lambda}$ has an associated RKHS, $\mathcal{H}_{k,\lambda}$ (endowed with its intrinsic norm $\|\cdot\|_{\mathcal{H}_{k,\lambda}}$), and the corresponding probability distance $d_{k,\lambda}$. For $\operatorname{P},\operatorname{Q}\in\mathcal{M}_{\operatorname{p}}(\mathcal{X})$, we want to test $\operatorname{H}_{0}:\operatorname{P}=\operatorname{Q}$ using the distances within the collection $\left\{d_{k,\lambda}:\lambda\in\Lambda\right\}$. The current theoretical framework does not support the automatic (data-driven) choice of $\lambda\in\Lambda$, since the asymptotic theory is mainly developed for a fixed kernel, corresponding to a specific value of $\lambda$. However, the choice of $\lambda$ is a non-trivial and sensitive issue with no obvious best solution, and which might affect the test performance. 

					There are various interesting proposals to deal with this problem in practice: the median heuristic of \cite{Borgwardt2006}; sample-splitting and optimization methods in \cite{Balakrishnan2012} and \cite{Gretton2020}; and aggregation methods such as \cite{Balakrishnan2012}. In this paper we explore an alternative to avoid making a parametric decision or splitting the data set. Our proposal can be included within the aggregative methods: we combine the information provided by different kernels by taking the supremum over the induced kernel metrics. Specifically, we use the metric
					that \enquote{best separates} $\operatorname{P}$ and $\operatorname{Q}$, that is, the supremum of all kernel distances given by
					\begin{equation}\label{definition.SKD}
						d_{k,\Lambda}(\operatorname{P},\operatorname{Q})=\underset{\lambda\in\Lambda}{\operatorname{sup}}\,\left(d_{k,\lambda}(\operatorname{P},\operatorname{Q})\right)=\underset{\lambda\in\Lambda}{\operatorname{sup}}\,\left(\left\|\mu_{\operatorname{P}}^{\lambda}-\mu_{\operatorname{Q}}^{\lambda}\right\|_{\mathcal{H}_{k,\lambda}}\right),\quad\operatorname{P},\ \operatorname{Q}\in\mathcal{M}_{\operatorname{p}}(\mathcal{X}),
					\end{equation}
					where, for $\lambda\in\Lambda$, $\mu_{\operatorname{P}}^{\lambda}$ and $\mu_{\operatorname{Q}}^{\lambda}$ are the mean embeddings of $\operatorname{P}$ and $\operatorname{Q}$, respectively, in $\mathcal{H}_{k,\lambda}$. We call the quantity in \eqref{definition.SKD} the \textit{supremum (or uniform) kernel distance} of $\left\{k_{\lambda}:\lambda\in\Lambda\right\}$. Also, the \textit{uniform kernel trick} refers to the overall idea of using \eqref{definition.SKD} to eliminate the parameter in kernel-based statistics. Observe that $d_{k}$ in \eqref{definition.kernel.distance} is a particular case of $d_{k,\Lambda}$ in \eqref{definition.SKD} when $\Lambda$ is a single-element set.
					Therefore, all the results in this work can be applied for usual kernel distances. In addition, in the family $\left\{k_{\lambda}:\lambda\in\Lambda\right\}$ we can include kernels from different parametric families, which would generate more robust test statistics that might work well under many types of alternatives.

					The supremum kernel distance \eqref{definition.SKD} entails several advantages and some mathematical challenges: First, the kernel selection problem is considerably simplified and solved in a natural way. Additionally, the approach is general enough to be applied in infinite-dimensional settings as FDA. This is interesting since in FDA there are only a few homogeneity tests in the literature. Some of them have been developed in the setting of ANOVA models (involving several samples) under homoscedasticity (equal covariance operators of the involved processes) and Gaussian assumptions. Hence, the current methodologies amount to testing the null hypothesis of equal means in all the populations; see, e.g., \cite{Cuevas2004} for an early contribution and \cite{Zhang2014} for a broader perspective. Our proposal is therefore quite related to more general approaches, not requiring any homoscedasticity assumption and still valid for a FDA framework. Examples of such similar tests are \cite{Ghosh2016} and \cite{Hall2007}, as well as the random projections-based methodology in \cite{Cuesta2007}.

					The inclusion of the supremum in \eqref{definition.SKD} represents an additional difficulty. The asymptotic properties of the test statistic based on \eqref{definition.SKD} are derived by following a different strategy from that of \cite{Borgwardt2006} and later works. The methodology proposed here allows us to cope with the supremum and applies directly to the case of unequal sample sizes. In short, our approach can be described as follows: First, we consider plug-in estimators of the kernel distances, obtained by replacing the unknown distributions by their empirical counterparts. Then, we use the powerful theory of empirical processes together with some recent results on the differentiability of the supremum (see \cite{Carcamo2020}) and functional Karhunen-Loève expansions of the underlying processes. These developments entail several technical difficulties from the mathematical point of view. However, they are worthwhile since they allow us to analyze the asymptotic behavior, under both the null and the alternative hypothesis, of the two-sample test based on \eqref{definition.SKD}.
				\subsubsection*{The organization of this paper}
					In Section \ref{Section_preliminaries} we provide some preliminaries regarding RKHS basics and empirical processes. While most of this background is well-known or can be found in the literature, it is included here to introduce the necessary notation and make the paper as self-contained as possible. Section \ref{Section_main_results} contains the main theoretical contributions. First, we obtain a Donsker property for unions of unit balls in RKHS that could be of independent interest. We establish the asymptotic validity under the null hypothesis of the two-sample test based on the distance \eqref{definition.SKD}. The asymptotic statistical power (i.e., the behaviour under the alternative hypothesis of non-homogeneity) is also analysed. An empirical study, comparing the uniform kernel test with some other competitors is presented in Section \ref{sec:empirical}. In the scenarios we have considered, SKD is competitive with the other kernel-based methods, especially in the case of heteroscedastic populations. However, given the limited nature of the study, we cannot conclude that our proposal unequivocally outperforms existing approaches. 
	\section{Preliminaries}\label{Section_preliminaries}
		In this section we describe various tools that we use throughout this work.
			\subsubsection*{Reproducing kernel Hilbert spaces (RKHS)}
				The theory of RKHS plays a relevant role in this paper. This is a classical and well-known topic; see \cite[Appendix F]{Janson1997} for a brief account of the RKHS theory and \cite{Berlinet&Thomas-Agnan2011} or \cite{Eubank&Hsing2015} for a statistical perspective. Hence, we only mention what is strictly necessary for later use.
				Let $\mathcal{X}$ be a topological space and $k:\mathcal{X}\times\mathcal{X}\rightarrow\mathbb{R}$ a \textit{kernel}, that is, a symmetric and positive semi-definite function. Let us consider $\mathcal{H}^{0}_{k}$, the pre-Hilbert space of all finite linear combinations $g(\cdot)=\sum_{i=1}^{n}\alpha_{i}\,k\left(x_{i},\cdot\right)$ (with $\alpha_{i}\in\mathbb{R}$, $n\in\mathbb{N}$ and $x_{i}\in\mathcal{X}$), endowed with the inner product
				\begin{equation}\label{eq:inner0}
					\left\langle\sum_{i=1}^{n}\alpha_{i}\,k\left(x_{i},\cdot\right),\sum_{j=1}^{m}\beta_{j}\,k\left(x_{j},\cdot\right)\right\rangle_{\mathcal H_{k}}=\sum_{i,j}\alpha_{i}\,\beta_{j}\,k\left(x_{i},x_{j}\right).
				\end{equation}
				The RKHS $\mathcal{H}_{k}$ is defined as the completion of $\mathcal{H}_{k}^{0}$; see \cite[Chapter 1]{Berlinet&Thomas-Agnan2011}. The inner product $\langle\cdot,\cdot\rangle_{k}$ in $\mathcal{H}_{k}$ is obtained through \eqref{eq:inner0} in such a way that bilinearity is preserved. A key property of RKHS is the so-called \textit{reproducing property}:
				\begin{equation}\label{eq:reproducing}
					\langle f,k(x,\cdot)\rangle_{\mathcal{H}_{k}}=f(x),\quad \operatorname{for\ all}\ f\in\mathcal{H}_{k},\ x\in \mathcal{X}.
				\end{equation}
			\subsubsection*{Kernel distances as integral probability metrics} 	
				Each $\operatorname{P}\in\mathcal{M}_{\operatorname{p}}(\mathcal{X})$ (Borel probability measure on $\mathcal{X}$), can be seen as a linear functional on $\mathcal{H}_{k}$ via the mapping
				\begin{equation}\label{P-linear-functional}
					f\in\mathcal{H}_{k}\mapsto\operatorname{P}(f)=\int_{\mathcal{X}}f\,\operatorname{dP},
				\end{equation}
				whenever $\mathcal{H}_{k}\subset\operatorname{L}^{1}(\operatorname{P})$ (set of integrable variables with respect to $\operatorname{P}$). This condition is also equivalent to saying that the function $x\mapsto k(x,\cdot)$ is Pettis integrable (with respect to $\operatorname{P}$) and to the existence of the mean embedding $\mu_{\operatorname{P}}$ in \eqref{eq:muP} as an element of $\mathcal{H}_{k}$ fulfilling Riesz representation condition
				\begin{equation}\label{eq:mean-embedding-property}
					\operatorname{}P(f)=\inp{f}{\mu_{\operatorname{P}}}{\mathcal{H}_{k}},\quad\operatorname{for}\ f\in\mathcal{H}_{k}.
				\end{equation}
				Sufficient conditions guaranteeing the injectivity of the mean embedding transformation can be found in  \cite{Fukumizu2011}. Note that in \eqref{P-linear-functional} (and what follows) we use the standard notation in empirical processes theory: $\operatorname{P}(f)$ (or simply $\operatorname{P}f$) stands for the mathematical expectation of $f$ with respect to $\operatorname{P}$.
				
				The existence of the mean embedding implies that the kernel distance in \eqref{definition.kernel.distance}, as well as the supremum kernel distance in \eqref{definition.SKD}, are well-defined.  Indeed, they are	\textit{integral probability metrics}; see \cite{Muller1997}. To see this, let us consider the unit ball of $\mathcal{H}_{k}$, that is,
				\begin{equation}\label{unit-ball}
					\mathcal{F}_{k}=\left\{f\in\mathcal{H}_{k},\ \|f\|_{\mathcal{H}_{k}}\leq1\right\}.
				\end{equation}
				We have that
				\begin{equation}\label{kernel.integral-metric}
					\begin{split}
						\left\|\mu_{\operatorname{P}}-\mu_{\operatorname{Q}}\right\|_{\mathcal{H}_{k}}&=\underset{f\in\mathcal{F}_{k}}{\operatorname{sup}}\left(\left\langle f,\mu_{\operatorname{P}}-\mu_{\operatorname{Q}}\right\rangle_{\mathcal{H}_{k}}\right)\overset{(a)}{=}\underset{f\in\mathcal{F}_{k}}{\operatorname{sup}}\left(\left\langle f,\int_{\mathcal{X}}k(\cdot,x)\,\operatorname{d}(\operatorname{P}-\operatorname{Q})(x)\right\rangle_{\mathcal{H}_{k}}\right)\\
						&\overset{(b)}{=}\underset{f\in\mathcal{F}_{k}}{\operatorname{sup}}\left(\int_{\mathcal{X}}\langle f,k(\cdot,x)\rangle_{\mathcal {H}_{k}}\,\operatorname{d}(\operatorname{P}-\operatorname{Q})(x)\right)\overset{(c)}{=}\underset{f\in\mathcal{F}_{k}}{\operatorname{sup}}(\operatorname{P}(f)-\operatorname{Q}(f)),
					\end{split}
				\end{equation}
				where $(a)$ follows from the definition of mean embedding \eqref{eq:muP}, $(b)$ from Pettis integrability, and $(c)$ from the reproducing property \eqref{eq:reproducing}; see also \cite[Lemma 4]{Borgwardt2012}. Thus, the kernel distance \eqref{definition.kernel.distance} is the integral probability metric generated by the class $\mathcal{F}_{k}$ in \eqref{unit-ball}.
				Therefore, the supremum kernel distance \eqref{definition.SKD} admits the alternative representation
				\begin{equation}\label{SKD-integral}
					d_{k,\Lambda}(\operatorname{P},\operatorname{Q})=\underset{f\in\mathcal{F}_{k,\Lambda}}{\operatorname{sup}}(\operatorname{P}(f)-\operatorname{Q}(f)) \quad\text{with}\quad \mathcal{F}_{k,\Lambda}=\bigcup_{\lambda\in \Lambda}  \mathcal{F}_{k,\lambda},
				\end{equation}	
				where $\mathcal{F}_{k,\lambda}$ is the unit ball in the RKHS space associated with $k_{\lambda}$. In other words, $d_{k,\Lambda}$ is the integral probability metric defined through the union of unit balls of the whole family of RKHS constructed with $\left\{k_{\lambda}:\lambda\in \Lambda\right\}$.

				From the characterizations as integral probability metrics in \eqref{kernel.integral-metric} and \eqref{SKD-integral}, we conclude that $d_{k}$ and $d_{k,\Lambda}$  satisfy the properties of a pseudo-metric (non-negativeness, symmetry, triangular property). However, to ensure the \textit{identifiability property} of a metric $d$ (i.e., $d(\operatorname{P},\operatorname{Q})=0$ if and only if $\operatorname{P}=\operatorname{Q}$) additional conditions are needed. It can be checked that when $\mathcal{X}=\mathbb{R}^{d}$, identifiability is satisfied for the usual kernels (such as the Gaussian kernel in \eqref{eq:gk}). However, when $\mathcal{X}$ is infinite-dimensional this type of results are more complicated; see \cite{Duncan&Wynne2022} for a deep study of this topic for the Gaussian kernel \eqref{eq:gk}. More details can also be found in \cite{Fukumizu2010} and \cite{Fukumizu2011}.
			\subsubsection*{Plug-in estimators, empirical processes and Donsker classes of functions}
				A simple and natural estimator of the supremum kernel distance \eqref{definition.SKD} can be obtained by applying \textit{the plug-in principle} in \eqref{SKD-integral}. Given two independent samples $X_{1},\ldots,X_{n}$ and $Y_{1},\ldots,Y_{m}$ from $\operatorname{P}$ and $\operatorname{Q}$, respectively, we replace the unknown underlying probability measures $\operatorname{P}$ and $\operatorname{Q}$ with the observed empirical counterparts,
				$$\mathbb{P}_{n}=\frac{1}{n}\,\sum_{i=1}^{n}\delta_{X_{i}},\quad\mathbb{Q}_{m}=\frac{1}{m}\,\sum_{i=1}^{m}\delta_{Y_{i}},$$
				$\delta_{a}$ being the unit point mass at $a$. This leads to the estimator of $d_{k,\Lambda}(\operatorname{P},\operatorname{Q})$ in \eqref{SKD-integral} given by
				\begin{equation}\label{plug-in_statistic}
					d_{k,\Lambda}\left(\mathbb{P}_{n},\mathbb{Q}_{m}\right)=\underset{f\in\mathcal{F}_{k,\Lambda}}{\operatorname{sup}}\,\left(\mathbb{P}_{n}(f)-\mathbb{Q}_{m}(f)\right)=\underset{f\in\mathcal{F}_{k,\Lambda}}{\operatorname{sup}}\,\left(\frac{1}{n}\,\sum_{i=1}^{n}f\left(X_{i}\right)-\frac{1}{m}\,\sum_{j=1}^{m}f\left(Y_{j}\right)\right).
				\end{equation}
				As a supremum over a class of functions is involved in \eqref{plug-in_statistic}, the theory of empirical processes comes into play naturally.  Given a collection of functions $\mathcal{F}$, we recall that the \textit{$\mathcal{F}$-indexed empirical process} (associated to $\operatorname{P}$) is $\mathbb{G}_{n}^{\operatorname{P}}=\sqrt{n}\,\left(\mathbb{P}_{n}-\operatorname{P}\right)$. The class $\mathcal{F}$ is called \textit{$\P$-Donsker} if $\mathbb{G}_{n}^{\operatorname{P}}\rightsquigarrow\mathbb{G}_{\operatorname{P}}$ in $\ell^{\infty}(\mathcal{F})$, the space of bounded real functionals defined on $\mathcal{F}$ with the supremum norm; see \cite{van_der_Vaart&Wellner2023}. Here, \enquote*{$\rightsquigarrow$} stands for weak convergence in $\ell^{\infty}(\mathcal{F})$ and $\mathbb{G}_{\operatorname{P}}$ is a \textit{$\operatorname{P}$-Brownian bridge}, that is, a zero-mean Gaussian process with covariance function
				\begin{equation*}
					\operatorname{E}\left[\mathbb{G}_{\operatorname{P}}\left(f_{1}\right)\,\mathbb{G}_{\operatorname{P}}\left(f_{2}\right)\right]=\operatorname{P}\left(f_{1}\,f_{2}\right)-\operatorname{P}\left(f_{1}\right)\,\operatorname{P}\left(f_{2}\right),\quad f_{1},f_{2}\in\mathcal{F}.
				\end{equation*}
				Additionally, $\mathcal{F}$ is \textit{universal Donsker} if it is $\operatorname{P}$-Donsker, for every $\operatorname{P}\in\mathcal{M}_{\operatorname{p}}(\mathcal{X})$.
	\section{Main results}\label{Section_main_results}
		In this section we first show that (unions of) unit balls of RKHS are universal Donsker under mild conditions. This is an important technical result of independent interest that is the starting point in the proofs of the asymptotic results. We analyze the asymptotic behaviour of the plug-in estimator \eqref{plug-in_statistic} of the supremum kernel distance in \eqref{definition.SKD} and \eqref{SKD-integral}. The results are quite general as $\operatorname{P}$ and $\operatorname{Q}$ are assumed to be Borel probability measures on a separable metric space. The proofs are based on empirical processes theory together with the (extended) delta method \cite[Theorem 2.1]{Shapiro1991} and some recent differentiability results for the supremum (\cite{Carcamo2020}). This \textit{differential approach} differs from previous methods (as those in \cite{Balakrishnan2012} or \cite{Borgwardt2012}) in which the theory of U-statistics is used to derive the asymptotic results. Our approach has some advantages: it is applicable to variables taking values in general spaces, including functional spaces, and the equal sample size constraint of previous works is removed. Furthermore, the results are applicable in other contexts (such as tests for equality between two copulas) by just changing the underlying stochastic process in the spirit of \cite{Carcamo2020}.
		
		Another essential difference between our methodology and other approaches is the way in which the tuning parameter $\lambda$ is treated. The asymptotic theory in \cite{Borgwardt2012} (and other related works) is derived for a fixed kernel, while the experiments incorporate the Gaussian kernel in \eqref{eq:gk} with a data-driven choice of $\lambda$. As pointed out by the authors, an automatic method for selecting $\lambda$ is an interesting area of research with some theoretical implications: setting the kernel using the sample being tested may cause changes to the asymptotic distribution. Regarding this, we note that our procedure to deal with the tuning parameter $\lambda$ is fully incorporated in the asymptotic analysis thanks to the use of the supremum kernel distance \eqref{definition.SKD}.
			\subsubsection*{The hypotheses}
				We list some assumptions for later reference. We briefly explain the meaning and implications of each of them. In what follows, $k$ is a kernel, $\left\{k_{\lambda}:\lambda\in\Lambda\right\}$ a family of kernels (which might come from different parametric families), and $\operatorname{P}$ and $\operatorname{Q}\in \mathcal{M}_{\operatorname{p}}(\mathcal{X})$, Borel probability measures defined on a space $\mathcal{X}$. In what follows we use the standard notation in functional analysis and operator theory; for $k_{1}$ and $k_{2}$ positive definite kernels on $\mathcal{X}$, we denote $k_{1}\ll k_{2}$ if and only if $k_{2}-k_{1}$ is a positive definite kernel; see \cite[Part I.7]{Aronszajn1950}.  
				\begin{description}
					\labitem{(Reg)}{itm:Reg} \textit{Regularity assumption.} $\mathcal{X}$ is a separable metric space and each kernel is continuous as a real function of one variable (with the other kept fixed).
					\labitem{(Dom)}{itm:Dom} \textit{Dominance assumption.} There exists a constant $c>0$ such that $k_{\lambda}\ll c\,k$, for all $\lambda\in\Lambda$. Further, $k$ is bounded on the diagonal, that is, $\underset{x\in\mathcal{X}}{\operatorname{sup}}\,(k(x,x))<\infty$.
					\labitem{(Ide)}{itm:Ide} \textit{Identifiability assumption.} If $\operatorname{P}\neq\operatorname{Q}$, there exists $\lambda\in\Lambda$ such that $\mu_{\operatorname{P}}^{\lambda}\neq\mu_{\operatorname{Q}}^{\lambda}$.
					\labitem{(Par)}{itm:Par} \textit{Continuous parametrization.} $\Lambda$ is a compact subset of $\mathbb{R}^{k}$ (with $k\in\mathbb{N}$) and, for a fixed $(x,y)\in\mathcal{X}\times\mathcal{X}$, the function $\lambda\mapsto k_{\lambda}(x,y)$ is continuous from $\Lambda$ to $\mathbb{R}$.
					\labitem{(Sam)}{itm:Sam} \textit{Sampling scheme.} The sampling scheme is balanced, that is, $\frac{n}{(n+m)}\to\theta$, with $\theta\in[0,1]$, as $n,m\to\infty$.
				\end{description}
				Assumptions \ref{itm:Reg} and \ref{itm:Dom}  together have important consequences. Firstly, they imply that $\mathcal{H}_{k,\lambda}$ is constituted by continuous and bounded functions (see \cite[Theorem 17]{Berlinet&Thomas-Agnan2011}), therefore measurable and integrable. Moreover, under these two conditions the mean embedding $\mu_{\operatorname{P}}^{\lambda}$ exits (for each $\operatorname{P}$ and $\lambda$). In particular, the supremum kernel distance \eqref{definition.SKD} is well-defined. \ref{itm:Reg} and \ref{itm:Dom} are also essential to show that the class  $\mathcal{F}_{k,\Lambda}$ in \eqref{SKD-integral} is universal Donsker, which is a key point in the proofs of the following theorems.

				Assumption \ref{itm:Ide} entails that $d_{k,\Lambda}(\operatorname{P},\operatorname{Q})>0$, whenever $\operatorname{P}\neq\operatorname{Q}$, i.e., the supremum kernel distance separates different probability measures. Therefore, $d_{k,\Lambda}$ in \eqref{definition.kernel.distance} is a proper metric on $\mathcal{M}_{\operatorname{p}}(\mathcal{X})$. Regarding this, we recall that a reproducing kernel $k$ is said to be \textit{characteristic} whenever $d_{k}(\operatorname{P},\operatorname{Q})=0$ if and only if $\operatorname{P}=\operatorname{Q}$, for all $\operatorname{P},\operatorname{Q}\in\mathcal{M}_{\operatorname{p}}(\mathcal{X})$; see \cite{Fukumizu2008}. This is equivalent to \enquote{integrally strictly positive definiteness},
				see \cite[Theorem 7]{Fukumizu2010}. Hence, \ref{itm:Ide} could be understood as \textit{the family $\left\{k_{\lambda}:\lambda\in\Lambda\right\}$ being characteristic} in the sense that for each pair of different measures there is a kernel in the family separating them. Observe that this condition is necessary to carry out the test $\operatorname{H}_{0}:\operatorname{P}=\operatorname{Q}$ by means of the statistic \eqref{plug-in_statistic}. Otherwise, the two-sample test that we propose only checks the weaker hypothesis $d_{k,\Lambda}(\operatorname{P},\operatorname{Q})=0$. Note that \ref{itm:Ide} is less demanding than asking for a specific kernel to be characteristic, a standard requirement on this topic. In infinite dimension, necessary and sufficient conditions for the Gaussian kernel to be characteristic are given in \cite{Duncan&Wynne2022}. We also observe that \ref{itm:Ide} is not specifically required to obtain the asymptotic distribution under $\operatorname{H}_{0}$ in Theorem \ref{Theorem_asymptotic_H0}.

				Finally, \ref{itm:Par} is a technical requirement to derive the asymptotic distribution of the test statistic under the alternative hypothesis using the results in \cite{Carcamo2020}. \ref{itm:Sam} is necessary for the combination of the associated empirical processes to converge.
			\subsubsection*{Examples of families of kernels}
				The hypotheses above can be verified for most families of kernels that are used in practice by properly choosing the parameter space. The most demanding assumption about the kernel family is perhaps \ref{itm:Dom}.
				In particular, this always ensures the applicability of the results, both in high and infinite dimensions, since the practical implementation of the procedure is carried out by choosing a grid of points in the parameter space; see Section \ref{sec:empirical}.
				
				When $\mathcal{X}=\mathbb{R}^{d}$, a finite-dimensional space, the usual parametric families of kernels often generate a nested collection of RKHS; see \cite{ZhangZhao2013}. This means that for $\lambda_{1},\lambda_{2}\in\Lambda$, there exists a constant $c=c\left(\lambda_{1},\lambda_{2},d\right)$ such that $k_{\lambda_{1}}\ll c\,k_{\lambda_{2}}$ (or the other way around). In such cases, \ref{itm:Dom} is valid for a compact subset $\Lambda$ of the whole parametric space by using one of the kernels of the family as the bounding kernel $k$ in \ref{itm:Dom}. Some important examples included in this setting are the families of Gaussian and Laplacian kernels, inverse multiquadrics kernels, B-spline kernels, Matérn kernels, among others; see \cite{Sriperumbudur2016} and \cite[Theorems 3.5, 3.6, and 3.7]{ZhangZhao2013}. 

				Nevertheless, the problem is more delicate in infinite dimension. If $\mathcal{X}=\mathbb{R}^{d}$ and for the usual parametric families of kernels, the best constant $c$ in \enquote{inequalities} of the form $k_{\lambda_{1}}\ll c\,k_{\lambda_{2}}$ depends on the dimension $d$ and blows up when $d$ goes to infinity; see \cite[Theorems 3.5 and 3.6]{ZhangZhao2013}. Therefore, when the domain is functional (for instance, if $\mathcal{X}=\operatorname{L}^{2}([0,1])$), the task of finding a dominating kernel is more involved. An example can be built when the parameter space
				$$\Lambda=\left\{m\in\operatorname{L}^{2}([0,1])\ \operatorname{:} m\ \operatorname{absolutely}\ \operatorname{continuous\ with}\ \int_{0}^{1}\left|m^{\prime}\right|^{2}<1\right\},$$
				is the unit ball of the Cameron-Martin space associated to the Wiener measure in $\mathcal{C}([0,1])$, the space of continuous functions on $[0,1]$. A family of kernels $k_{m}(x,y)$ fulfilling \ref{itm:Dom} and \ref{itm:Ide} can be constructed using Minlos-Sazanov Theorem (see e.g., \cite[Th. 24]{Duncan&Wynne2022}) and relying on ideas by \cite[Chapter 10]{Wendland2004} and \cite[Prop. 3.1]{ZhangZhao2013}.
				
				Additionally, we observe that \ref{itm:Dom} is fulfilled for families of positive linear (or convex) combinations of a finite family of kernels. In this example, the set of parameter $\Lambda$ is given by the weights of the considered combinations; see \cite{Balakrishnan2012}. We finally refer to \cite[Chapter 7]{Berlinet&Thomas-Agnan2011} and \cite[Chapter 4]{Paulsen2016} for a wider catalog of families of kernels within this context.
			\subsubsection*{A Donsker property for units balls in RKHS}
				Establishing that a class of functions is (uniform) Donsker has important consequences. This property is equivalent to having an empirical central limit theorem, which is at the heart of most asymptotic results in statistics. Therefore, this kind of Donsker-type results are relevant by themselves and of independent interest. For example, in \cite[Theorem 4.3]{Sriperumbudur2016} (see also \cite{Gine&Nickl2008}, \cite{Gine&Nickl2016}) it is shown that $\mathcal{F}_{k,\Lambda}$ in \eqref{SKD-integral} is Donsker for some specific finite-dimensional parametric families and for a suitable subset of $\Lambda$. Then, this result is applied to derive asymptotic distributions of kernel density estimators. In \cite{Sriperumbudur2016}, the proofs of the Donsker property for RKHS unit balls are obtained when $\mathcal{X}=\mathbb{R}^{d}$ by direct covering (entropy-based) arguments. The underlying bounds in these references depend on the dimension $d$. Therefore, it seems difficult to extend these Donsker-type statements to the infinite-dimensional case. However, Theorem \ref{Theorem:Donsker} below is suitable for the general framework where $\mathcal{X}$ might be an infinite-dimensional space, and thus useful in statistical problems with functional data.
				
				The following theorem establishes that unit balls (and even the union of units balls) of RKHS are universal Donsker. In the first part of the proof (in Section  \ref{sec:Proofs}) we use \cite[Theorem 1.1]{Marcus1985}, while in the second one we show that the union of unit balls is included in a ball of the space $\mathcal{H}_k$ by using  Aronszajn's inclusion theorem (\cite[Theorem I]{Aronszajn1950}). 
				\begin{Th} \label{Theorem:Donsker}
					Assume that the kernel $k$ satisfies \ref{itm:Reg} and it is bounded on the diagonal, that is, $\underset{x\in\mathcal{X}}{\operatorname{sup}}\,(k(x,x))<\infty$.
					Then, the class $\mathcal{F}_{k}$ in \eqref{unit-ball} (i.e., the unit ball of $\mathcal{H}_k$) is universal Donsker. Moreover, if $\left\{k_{\lambda}:\lambda\in\Lambda\right\}$ and $k$ satisfy \ref{itm:Dom}, then the union $\mathcal{F}_{k,\Lambda}$ in \eqref{SKD-integral} is universal Donsker as well.
					%
				\end{Th}
				This theorem extends \cite[Theorem 4.3]{Sriperumbudur2016}, where the Donsker property was shown under more demanding analytical conditions, to any family of kernels satisfying \ref{itm:Dom}.
			\subsubsection*{Asymptotic behaviour under the null hypothesis, $\operatorname{P}=\operatorname{Q}$}
				The next theorem provides the asymptotic distribution of the (normalized) estimator of the supremum kernel distance \eqref{definition.SKD} when the two samples come from the same distribution. In the statement of the following results, $\mathbb{G}_{\operatorname{P}}$ and $\mathbb{G}_{\operatorname{Q}}$ are $\mathcal{F}_{k,\Lambda}$-indexed $\operatorname{P}$ and $\operatorname{Q}$ Brownian bridges, respectively (see Section \ref{Section_preliminaries}), \enquote*{$\rightsquigarrow$} stands for the usual convergence in distribution of (real) random variables, and $\mathcal{H}_{k,\lambda}^{\ast}$ denotes the dual space of $\mathcal{H}_{k,\lambda}$.
				\begin{Th}\label{Theorem_asymptotic_H0}
					Let us assume that \ref{itm:Reg}, \ref{itm:Dom} and \ref{itm:Sam} hold. If $\operatorname{P}=\operatorname{Q}$, the statistic \eqref{plug-in_statistic} satisfies that
					\begin{equation}\label{eq:theorem-H0}
						\sqrt{\frac{n\,m}{n+m}}\,d_{k,\Lambda}\left(\mathbb{P}_{n},\mathbb{Q}_{m}\right)\rightsquigarrow\underset{\lambda\in\Lambda}{\operatorname{sup}}\,\left(\left(\sum_{j\in\mathbb{N}}Z_{j,\lambda}^{2}\right)^{\rfrac{1}{2}}\right), \quad n,m\to\infty,
					\end{equation}
					where $d_{k,\Lambda}$ is defined in \eqref{SKD-integral}, $Z_{j,\lambda}=\left\langle\mathbb{G}_{\operatorname{P}},\varphi_{j,\lambda}\right\rangle_{\mathcal{H}_{k,\lambda}^{\ast}}$ (for each $\lambda\in\Lambda$ and $j\in\mathbb{N}$) and $\varphi_{j,\lambda}$ is the $j$-th eigenfunction of the covariance operator of $\mathbb{G}_{\operatorname{P}}$ on $\mathcal{H}_{k,\lambda}^{\ast}$.
	
					Moreover, $\left\{Z_{j,\lambda}\right\}_{j\in\mathbb{N},\lambda\in\Lambda}$ are jointly Gaussian and for a fixed $\lambda\in\Lambda$, $\left\{Z_{j,\lambda}\right\}_{j\in\mathbb{N}}$ are independent with $Z_{j,\lambda}\sim\mathcal{N}\left(0,\beta_{j,\lambda}\right)$, where $\beta_{j,\lambda}$ is the eigenvalue associated to $\varphi_{j,\lambda}$.
				\end{Th}
				In the first step of the proof of this theorem we use Theorem~\ref{Theorem:Donsker} to derive the weak convergence of the underlying process. The rest of the proof is rather technical. The basic ideas are as follows: we use of the continuous mapping theorem to obtain the convergence of the statistic; subsequently, we apply a functional Karhunen-Loève-type theorem in the dual space $\mathcal{H}_{k,\lambda}^{\ast}$ (Lemma \ref{Lemma.KL} in Section \ref{sec:Proofs}) to the resulting limiting process to achieve \eqref{eq:theorem-H0}. Note that in the family $\left\{k_{\lambda}:\lambda\in\Lambda\right\}$ we can include kernels from different parametric families or mixtures of kernels from distinct families in order to \textit{robustify} the test statistic.
	
				Theorem \ref{Theorem_asymptotic_H0} complements in several directions other previous works on this topic, starting from \cite[Th. 8]{Borgwardt2006}. See also \cite{Guo2022}, \cite{Smaga&Zhang2022} for more recent references.
			\subsubsection*{Asymptotic behaviour under the alternative, $\operatorname{P}\neq\operatorname{Q}$}
				The following theorem establishes the asymptotic distribution of (the normalized version) of \eqref{plug-in_statistic} under the alternative hypothesis of the homogeneity test. Therefore, it provides the consistency of the testing procedure based on the supremum kernel distance. Additionally, this result might be potentially useful in order to develop tests of \textit{almost homogeneity}, that is, problems in which we are interested in testing $\operatorname{H}_{0}:d_{k,\Lambda}(\operatorname{P},\operatorname{Q})\leq\varepsilon$ versus $\text{H}_{1}:d_{k,\Lambda}(\operatorname{P},\operatorname{Q})>\varepsilon$, for some $\varepsilon>0$. Analogously, this idea is also applicable to provide statistical evidence in favour of almost homogeneity when $\operatorname{H}_{0}$ and $\operatorname{H}_{1}$ above are interchanged. Related ideas can be found in \cite{Barrio-et-al-2020} and \cite{Dette&Kokot2022}.
				\begin{Th}\label{Theorem_asymptotic_H1}
					Let us assume that \ref{itm:Reg}, \ref{itm:Dom}, \ref{itm:Par}, \ref{itm:Ide} and \ref{itm:Sam} hold. If $\operatorname{P}\neq\operatorname{Q}$, we have that
					\begin{equation}\label{eq:theorem2}
						\sqrt{\frac{n\,m}{n+m}}\,\left(d_{k,\Lambda}\left(\mathbb{P}_{n},\mathbb{Q}_{m}\right)-d_{k,\Lambda}(\operatorname{P},\operatorname{Q})\right)\rightsquigarrow\underset{\lambda\in\Lambda_{0}}{\operatorname{sup}}\,\left(\mathbb{G}\left(h^{+,\lambda}\right)\right)=\underset{L}{\operatorname{sup}}\,(\mathbb{G}),
					\end{equation}
					where
					\begin{equation}\label{Witness}
						\mathbb{G}=\sqrt{1-\theta}\,\mathbb{G}_{\operatorname{P}}-\sqrt{\theta}\,\mathbb{G}_{\operatorname{Q}},\quad h^{+,\lambda}=\frac{\mu_{\operatorname{P}}^{\lambda}-\mu_{\operatorname{Q}}^{\lambda}}{\left\|\mu_{\operatorname{P}}^{\lambda}-\mu_{\operatorname{Q}}^{\lambda}\right\|_{\mathcal{H}_{k,\lambda}}},
					\end{equation}
					\begin{equation}\label{Witness-sets}
						\Lambda_{0}=\left\{\lambda\in\Lambda: \left\|\mu_{\operatorname{P}}^{\lambda}-\mu_{\operatorname{Q}}^{\lambda}\right\|_{\mathcal{H}_{k,\lambda}}=d_{k,\Lambda}(\operatorname{P},\operatorname{Q})\right\}\quad\operatorname{and}\quad L=\left\{h^{+,\lambda}:\lambda\in\Lambda_{0}\right\}.
					\end{equation}
				\end{Th}
				Theorem \ref{Theorem_asymptotic_H1} directly provides the consistency of the homogeneity test based on the supremum kernel distance $d_{k,\Lambda}$ in \eqref{definition.SKD}. We also observe that $\mathbb{G}$ is a zero mean Gaussian process indexed by $\mathcal{F}_{k,\Lambda}$. Further, $h^{+,\lambda}$ is called \textit{witness function} in \cite{Borgwardt2006} as the maximum mean discrepancy over $\mathcal{F}_{k,\lambda}$ is attained at this element, that is, $\operatorname{P}\left(h^{+,\lambda}\right)-\operatorname{Q}\left(h^{+,\lambda}\right)=\|\mu_{\operatorname{P}}^{\lambda}-\mu_{\operatorname{Q}}^{\lambda}\|_{\mathcal{H}_{k,\lambda}}$. Therefore, the limit in \eqref{eq:theorem2} corresponds to the supremum of $\mathbb{G}$ over the set of witness functions for which the value of the uniform kernel distance is achieved. Regarding the proof of Theorem \ref{Theorem_asymptotic_H1}, we mention that the extended delta method (see \cite[Theorem 2.1]{Shapiro1991}) plays a key role. First, we use Theorem \ref{Theorem:Donsker} to show that $\mathbb{G}$ is the limit of the underlying process. Later, we adapt some ideas from \cite{Carcamo2020} to derive \eqref{eq:theorem2}.
				
				The following result is a direct consequence of Theorem \ref{Theorem_asymptotic_H1} when the family of kernels has a single element, $k$.
				\begin{Cor}\label{Corollary-H1}
					Let us assume that \ref{itm:Reg}, \ref{itm:Dom} and \ref{itm:Sam} hold. Further, we assume that $k$ is characteristic. If $\operatorname{P}\neq\operatorname{Q}$, we have that
					\begin{equation}\label{eq:corolary1}
						\sqrt{\frac{n\,m}{n+m}}\,\left(d_{k}\left(\mathbb{P}_{n},\mathbb{Q}_{m}\right)-d_{k}(\operatorname{P},\operatorname{Q})\right)\rightsquigarrow\mathbb{G}\left(h^{+}\right),
					\end{equation}
					where $\mathbb{G}$ is in \eqref{Witness} and
					\begin{equation}\label{h+funtion}
						h^{+}=\frac{\mu_{\operatorname{P}}-\mu_{\operatorname{Q}}}{\left\|\mu_{\operatorname{P}}-\mu_{\operatorname{Q}}\right\|_{\mathcal{H}_{k}}}.
					\end{equation}
					In particular, the distribution of $\mathbb{G}\left(h^{+}\right)$ is normal with mean zero and variance $\Var\left(\mathbb{G}\left(h^{+}\right)\right)=(1-\theta)\,\Var_{\operatorname{P}}\left(h^{+}\right)+\theta\,\Var_{\operatorname{Q}}\left(h^{+}\right)$.
				\end{Cor}
				Corollary \ref{Corollary-H1} extends some previous results in which it is assumed that $n=m$; see \cite[Th. 8]{Borgwardt2006}, \cite[Th. 2.5]{Borgwardt2006a}, and \cite[Th. 16]{Duncan&Wynne2022}.
	\section{Empirical results}\label{sec:empirical}
		The aim of this section is to provide some insight about the performance of the two-sample test based on the SKD in \eqref{definition.SKD}, both from simulations and real world data sets.
			\subsubsection*{The purpose of these experiments and the methods under study}
				In the same spirit as \cite{Borgwardt2012} or \cite[Section 8.1]{Fukumizu2007}, we emphasize the interest of a new homogeneity test (based on kernel distances), suitable for high-dimensional data and not suffering from the degradation of classical two-sample tests when the dimension increases. Additionally, we show the advantages of avoiding the choice of the parameters in kernel distances, via our SKD proposal. The general idea is to check the SKD methodology as an attempt to robustify the test statistic against bad choices of the kernel or its parameter(s). In this empirical study we compare the following methods:
				\begin{description}
					\item [--] SKD: test based on the SKD in \eqref{definition.SKD} with a Gaussian kernel in \eqref{eq:gk}.
					\item [--] GKD: the kernel distance-based test of \cite{Fukumizu2007} with a data-driven choice of $\lambda$ in \eqref{eq:gk} as the median distance between points in the aggregate sample.
					\item [--] GKDSplit: test based on a Gaussian kernel distance where the estimation of the parameters is done by a splitting method; see \cite{Fukumizu2009}. The sample is divided into training and test subsamples to avoid data influence on the parameter selection.
					\item [--] GKDSplitOpt: test based on a kernel distance where the parameter estimation is detailed in \cite[Section 3]{Balakrishnan2012}. The data is divided into training and test sets. The target parameters are the coefficients of a convex combination of a finite family of kernels. The weights of the combination are selected to maximize the ratio between the empirical bias-corrected kernel distance and the standard deviation of the associated asymptotic distribution of the centered, normalized empirical distance.
					\item [--] ET: the energy test, a popular choice in this type of problems; see \cite{Rizzo&Szekely2017}, \cite{EnergyPackage2022}.
				\end{description}
				In our view, the \enquote{data-splitting} proposals are
				based on natural ideas that deserve attention. Still, there are some open issues to clarify, especially regarding the optimal splitting of the sample and the asymptotic behaviour of the resulting data-driven tests. We hope that this study could encourage further research along these lines. As for the energy test ET (\cite{Rizzo&Szekely2017}), we have included it in the study because it is based in a successful statistical methodology, ultimately grounded on the underlying \enquote{distance covariance} association measure; in fact this method has become quite popular in  high-dimensional two sample problems, via the \verb!energy! package in \verb!R!. Following a suggestion from one of the reviewers, we have considered as well a variant of this method, which is based on a different distance between the sample points. The standard statistic in ET is calculated in terms of the Euclidean distance, which can be seen, by a duality reasoning explained in \cite{Rizzo&Szekely2017}, as a distance associated with the Brownian covariance. This suggests the possibility of broadening the choice of the distance to the whole range of the fractional Brownian motion with Hurst parameter $H\in(0,1)$; recall that $H={1}/{2}$ for the standard Brownian motion. In our case, the choice $H={3}/{4}$ led to results almost identical to those of $H={1}/{2}$. Perhaps if heavy-tailed distributions were involved, this new alternative could make a real difference.

				The present empirical study is intended as an illustration of our proposal. Therefore, it is far from exhaustive. In particular, in which concerns SKD method, there is a considerable room to check the influence of different grids for the optimization on $\lambda$. A much more detailed experiment, including additional models and competitors might be worthwhile, but this is beyond the scope of this paper.
			\subsubsection*{The models}
				We include the models in \cite{Fukumizu2007}, based on Gaussian distributions in high dimension with different means and diagonal covariance matrices. In addition, we also consider a new scenario with functional data corresponding to trajectories of Gaussian processes in $\operatorname{L}^{2}([0,1])$. We note that all the considered tests can be applied in the functional setting: GKD, SKD, GKDSplit and GKDSplitOpt are based on the aggregated matrix $\left(k_{\lambda}\left(Z_{i},Z_{j}\right)\right)_{i,j=1}^{n+m}$, where $Z_{l}=X_{l}$ for $l\in\{1,\ldots,n\}$ and $Z_{n+l}=Y_{l}$ for $l\in\{1,\ldots,m\}$. ET uses cross-distances between the data in the sample space.
				
				More specifically, our simulation experiments are grouped in three blocks, respectively corresponding to different versions of homoscedasticity (Experiment 1) and heteroscedasticity (Experiment 2), plus a functional real data example.
				\begin{description}
					\item[\textbf{Experiment 1.}]\textsl{Different means, homoscedastic case}
					\begin{description}
						\item[\textsl{Model 1.1}] \textit{White noise.} We consider $\operatorname{P}\sim\mathcal{N}(0,\mathbb{I})$ and $\operatorname{Q}\sim\mathcal{N}(\mu\,\mathbf{1},\mathbb{I})$, where $\mathbf{1}=\left(1,\overset{d)}{\ldots},1\right)^{\top}$ (the superindex denotes the transpose) and $\mathbb{I}$ is the $d\times d$ identity matrix. In this model we deal with two multivariate Gaussian distributions with identity covariance in large dimension: $\operatorname{P}$ is standard and $\operatorname{Q}$ has mean $\mu\,\mathbf{1}$. Hence, $\operatorname{Q}$ is a shifted version of $\operatorname{P}$ translated $\frac{\mu}{\sqrt{d}}$ units in the direction given by the vector $\mathbf{1}$. The parameter $\mu$ takes the values $0$ (null hypothesis), $0.01$, $0.02$ and $0.05$ (alternative hypothesis).
						\item[\textsl{Model 1.2}] \textit{Functional data.} In this case $\operatorname{P}\sim\mathcal{G}(0,\gamma)$ and $\operatorname{Q}\sim\mathcal{G}(\mu\,\mathbf{1},\gamma)$, where $\mathcal{G}$ stands for a Gaussian process in $\operatorname{L}^{2}([0,1])$. The first parameter is the mean function and the second the covariance function. Here, $\mathbf{1}$ is the function identically equal to 1 and $\gamma\left(t_{1},t_{2}\right)=\operatorname{exp}\left(-0.5\,\left|t_{1}-t_{2}\right|\right)$. In this model, the \enquote{dimension} refers to the size of the grid used to approximate the process. The parameter $\mu$ takes the values $0$ (null hypothesis), $0.01$, $0.05$ and $0.2$ (alternative hypothesis).
					\end{description}
					\item[\textbf{Experiment 2.}]\textsl{Equal means, heteroscedastic cases}
					\begin{description}
						\item[\textsl{Model 2.1}] \textit{Spread white noise.} We consider $\operatorname{P}\sim\mathcal{N}(0,\mathbb{I})$ and $\operatorname{Q}\sim\mathcal{N}(0,\sigma^{2}\,\mathbb{I})$. The measure $\operatorname{P}$ corresponds to a standard multidimensional Gaussian distribution and $\operatorname{Q}$ to $\sigma$ times $\operatorname{P}$. The parameter $\sigma^{2}$ takes the values $10^{0.01}$ and $10^{0.02}$. This scenario introduces different alternative hypotheses from those in Model 1.1. In this example, $\operatorname{P}$ is more concentrated around the mean than $\operatorname{Q}$.
						\item[\textsl{Model 2.2}] \textit{Equicorrelated marginals.} Here, $\operatorname{P}\sim\mathcal{N}(0,\mathbb{I})$ and $\operatorname{Q}\sim\mathcal{N}(0,\Sigma)$, where $\Sigma=\rho\,\left(\mathbf{1}\,\mathbf{1}^{\operatorname{T}}-\mathbb{I}\right)+\mathbb{I}$, with $\rho\in\{0.005,0.01, 0.02, 0.05\}$. This scenario includes another different alternative from the ones in Model 1.1. In this case, the difference between $\operatorname{P}$ and $\operatorname{Q}$ lies on the (linear) dependence structure of the marginals.
					\end{description}
				\end{description}
			\subsubsection*{Some  technical aspects}
				Throughout this study we restrict ourselves to the family of Gaussian kernels in \eqref{eq:gk}, where $\mathcal{X}=\mathbb{R}^{d}$ or $\operatorname{L}^{2}([0,1])$. In this case, it is easy to show that the kernel distance $d_{k,\lambda}(\operatorname{P},\operatorname{Q})\to0$,  when $\lambda\to0$ or $\lambda\to\infty$ (and discrete part of $\operatorname{P}$ and $\operatorname{Q}$ is null). Given two random samples $X_{1},\ldots,X_{n}\sim\operatorname{P}$ and $Y_{1},\ldots,Y_{m}\sim\operatorname{Q}$, due to the fact that the empirical measures $\mathbb{P}_{n}$ and $\mathbb{Q}_{m}$ are discrete, the kernel distance
				$$d_{k,\lambda}\left(\mathbb{P}_{n},\mathbb{Q}_{m}\right)=O\left(\sqrt{\frac{1}{n}+\frac{1}{m}}\right),\quad\text{operatorname}\ \lambda\to\infty.$$
				This means that for small sample sizes, the plug-in estimator of the distance does not properly approximate its population counterpart. In particular, the maximum of the empirical distance is usually attained \enquote{at the tail}, i.e., on the extremes of the target interval for $\lambda$. This drawback is inherent to the classical kernel distance although it has not been explicitly mentioned in the literature. Therefore, it is convenient to slightly modify the kernel to make the empirical distance behave the same as in the continuous case for small sample sizes.  We propose to use a smoothed Gaussian kernel for this experiments given by
				$$k_{\lambda}(x,y)=\operatorname{exp}\left(-\lambda\,\left(\|x-y\|^{2}+0.1\,\left(\|x\|^{2}+\|y\|^{2}\right)\right)\right).$$
				This regularization is common in harmonic analysis to approximate the Dirac delta in spaces of distributions via smooth functions, called \textit{mollifiers}. It could be seen as an ad-hoc correction to improve the approximation of the maximum of the estimated kernel distance to the corresponding \enquote{true} population maximum. The smoothing process can be eliminated when sample sizes are sufficiently large.
				
				As shown in the literature, a data-driven choice of $\lambda$ seems to have a good practical behaviour. Specifically, in \cite{Fukumizu2007} (and in subsequent works), the value of $\lambda$ is the median distance between points in the aggregate sample. A theoretical consequence of this choice is that the asymptotic theory, derived in \cite{Fukumizu2007} under the assumption that $\lambda$ is fixed, does not longer apply to the data-driven case. Still, we include in our experiments, for comparison purposes, this data-driven choice as it is a common practice in the earlier literature. Since, to the best of our knowledge, the asymptotic distribution of the data-driven statistic is not known, we use a permutation test based on this statistic to obtain rejection regions rather than the other methods (Pearson curves, Gamma curves and bootstrap for U-statistics) explained in \cite{Fukumizu2007}.

				It is worth mentioning that in both tests  SKD and ET a permutation procedure has been used to approximate the corresponding distributions. This is the methodology used in the energy \verb!R!-package for the ET test and we have followed here the same strategy for the SKD test. Let us recall that our theoretical results provide the asymptotic distribution of this test, as well as its consistency (see Theorems \ref{Theorem_asymptotic_H0} and \ref{Theorem_asymptotic_H1}). However, the estimation of the quantiles of the limit distribution in \eqref{eq:theorem-H0} is far from trivial. As an additional complication, standard bootstrap approximation fails, as a consequence of the results in \cite{Fang&Santos2019}. This is why the permutation method appears as a natural choice. Further, the validity of the permutation test comes from the permutation invariance of the test statistic (see \cite{Lehmann&Romano2022}).

				Let us recall that the idea behind the SKD test is to dodge the parameter selection problem by considering \enquote{the whole parameter space}. Ideally, for the Gaussian kernel, an interval of the form $(0,\infty)$ could be considered in the SKD. As mentioned before, the extremes ($0$ and $\infty$) are not useful since the distance tends to zero when the parameter approaches to these values. Therefore, we use a parameter space of the form $\Lambda=[a, b]$, with $0<a<b<\infty$. It is important to note that here the values $a$ and $b$ cannot be properly considered as tuning parameters, since the test is not particularly sensitive to their choice, provided that the interval is large enough. As we have experimentally verified, $\Lambda=\left[10^{-4},0.1\right]$ is adequate to carry out the test. In practice, we employ a grid of $11$ points logarithmically separated between $10^{-4}$ and $0.1$. The goal is to approximate the value of the supremum by the maximum over finite subsets. The simulation outputs below are based on averages over $200$ replications. The permutation tests for GKD and SKD correspond to $B=5000$ permutations. As for ET, we use the function \verb!eqdist.etest! of the \verb!R!-package \cite{EnergyPackage2022}. Sample sizes are $n=m=250$ in all experiments. The effect of increasing the dimension $d$ is checked in the rank $d\in\{205,405,603,803,1003,1203,1401,$ $1601,1801,2001\}$. In all cases, the significance level of the test is set at $\alpha=0.05$.
			\subsubsection*{Outputs}
				\begin{figure}[!h]
					\centering
					\includegraphics[width=0.9\linewidth]{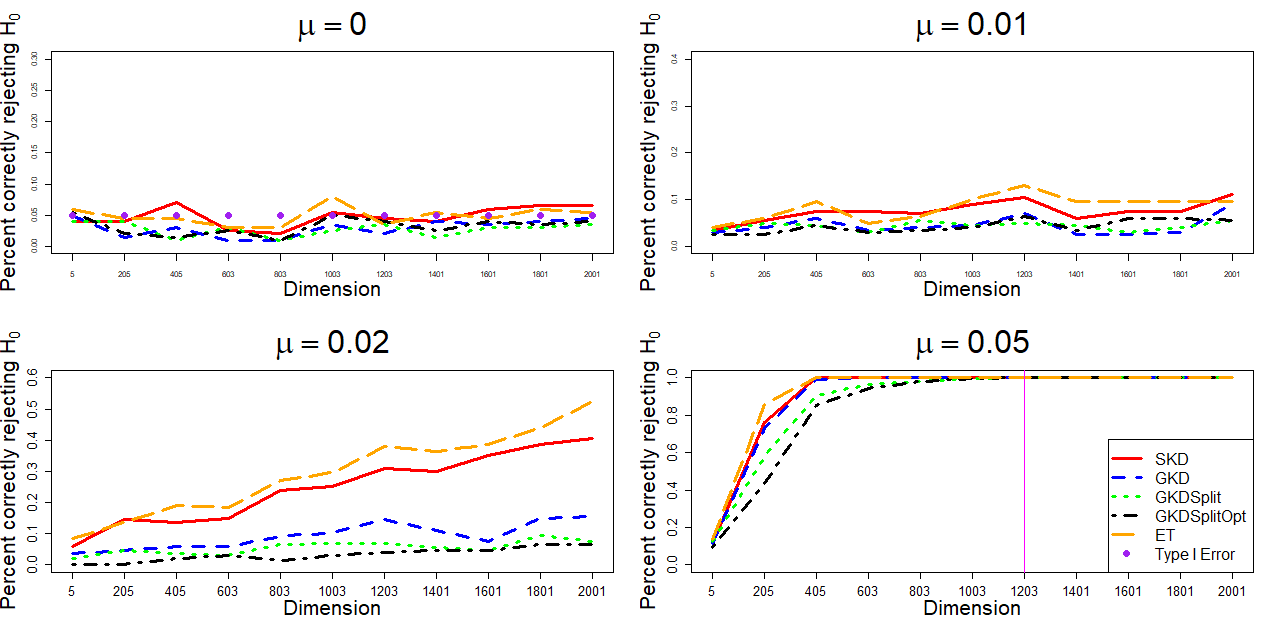}
					\caption{Performance of the tests under Model 1.1 with $\alpha=0.05$. Four values of $\mu$ are shown: $0$ (null hypothesis), $0.01$, $0.02$, and $0.05$ (alternative hypothesis).}
					\label{fig:Model1.1}
				\end{figure}
				Outputs from Model 1.1 are displayed in Fig. \ref{fig:Model1.1}. Tests calibration, i.e., the behaviour of the different tests under $\text{H}_0$, corresponds to the case $\mu=0$. We observe that the size of the test is reasonably well controlled by the five tests. Under the alternative hypothesis $\mu=0.01$ power curves oscillate slightly and are low. 
				
				\begin{figure}[!h]
					\centering
					\includegraphics[width=0.9\linewidth]{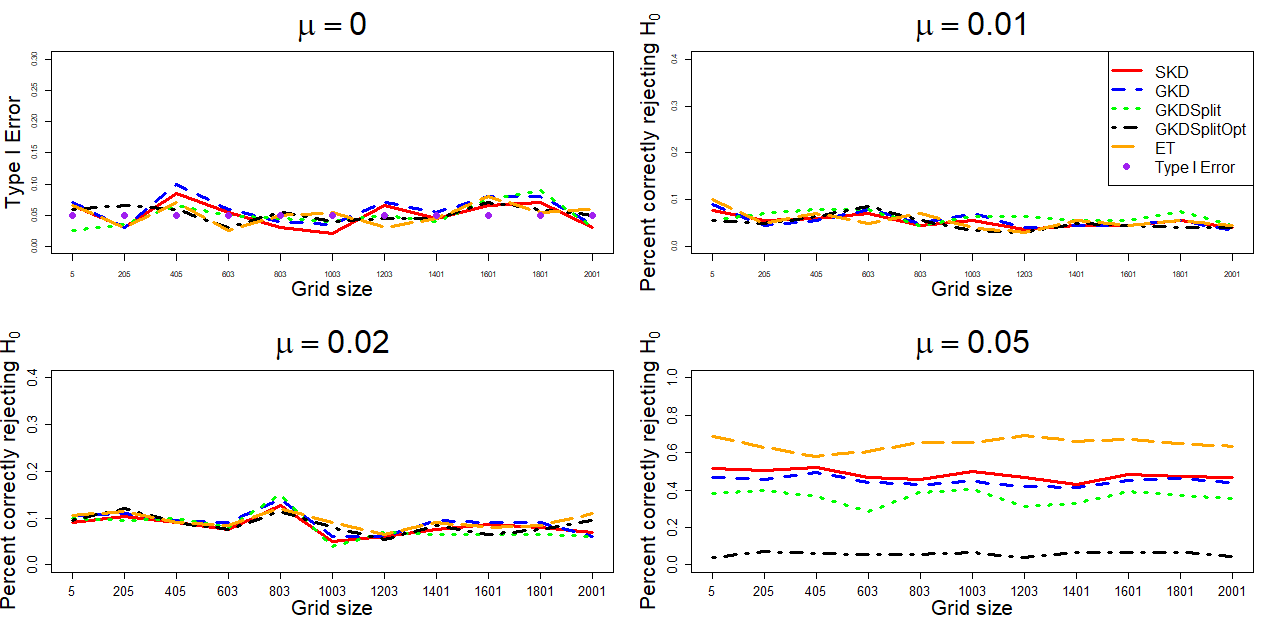}
					\caption{Performance of the tests under Model 1.2 with $\alpha=0.05$. Four values of $\mu$ are shown: $0$ (null hypothesis), $0.01$, $0.05$, and $0.2$ (alternative hypothesis).}
					\label{fig:Model1.2}
				\end{figure}
				Results of Model 1.2 are summarized in Fig. \ref{fig:Model1.2}. Test calibration outputs are depicted for $\mu=0$. As in the previous example, the second graph (case $\mu=0.01$) shows a relatively small power in all cases, since both distributions are very close to each other. A gain in power is observed for $\mu=0.05, 0.2$. The functional nature of the data is apparent in the fact that there is no clear pattern of \enquote{dimensionality blessing} associated with the increase of grid size. Indeed, unlike the other examples we are considering, the use of higher dimensional observations (a denser grid) does not entail a true gain in information, as grid observations are highly correlated, due to the continuity of the trajectories. ET obtains the best results under this scenario and SKD is competitive with the other kernel-based tests. Finally, it is surprising the absence of power of GKDSplitOpt in this model.
				
				\begin{figure}[!h]
					\centering
					\includegraphics[width=0.9\linewidth]{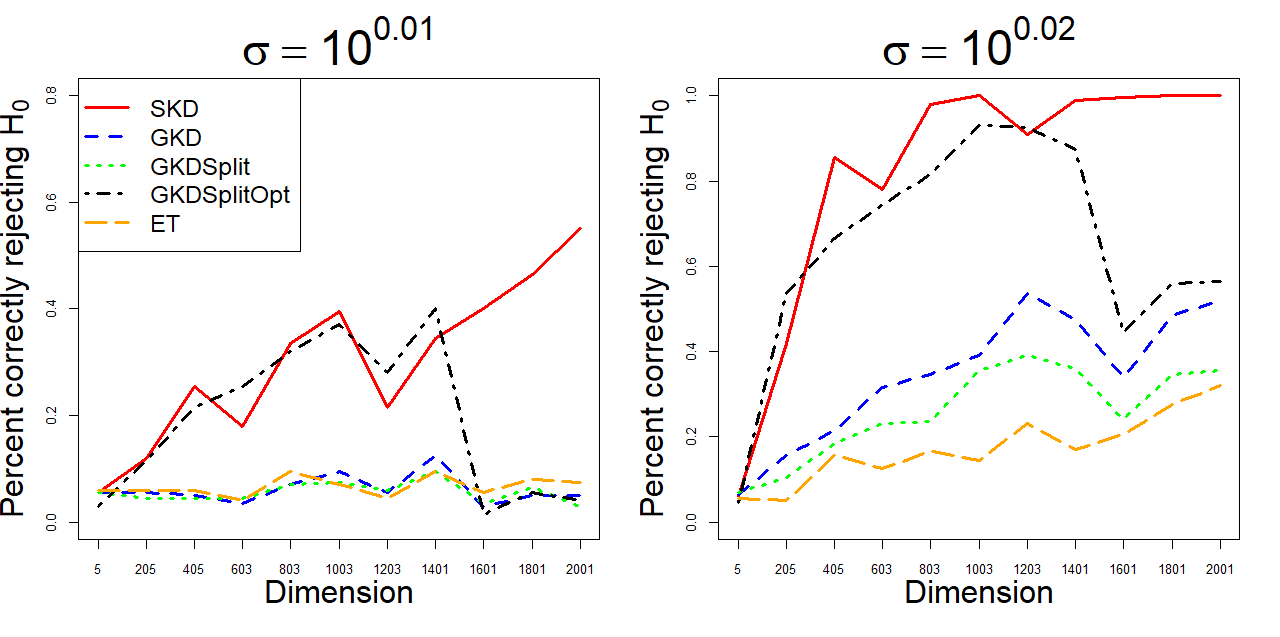}
					\caption{Performance of tests under Model 2.1 with $\alpha=0.05$. Two values of $\sigma^{2}$ are shown: $10^{0.01}$ and $10^{0.02}$ (alternative hypothesis).}
					\label{fig:Model2.1}
				\end{figure}
				The outputs from the heteroscedastic Model 2.1 are placed in Fig. \ref{fig:Model2.1}. Here GKDSplitOpt and SKD behave very well and clearly outperform the other methods. A plausible explanation for this difference is that the (data-driven) median-based selection of $\lambda$ of GKD is not a good choice for the heteroscedastic case when the value of the location parameter is the same in both populations. This heteroscedastic, same-location, scenario is also not the most favorable for the ET. Finally, it is noteworthy the loss of power of GKDSplitOpt from dimension 1601 onward.  Again, this might be due to the small sample sizes in relation to the dimension of the problem.
				
				\begin{figure}[!h]
					\centering
					\includegraphics[width=0.9\linewidth]{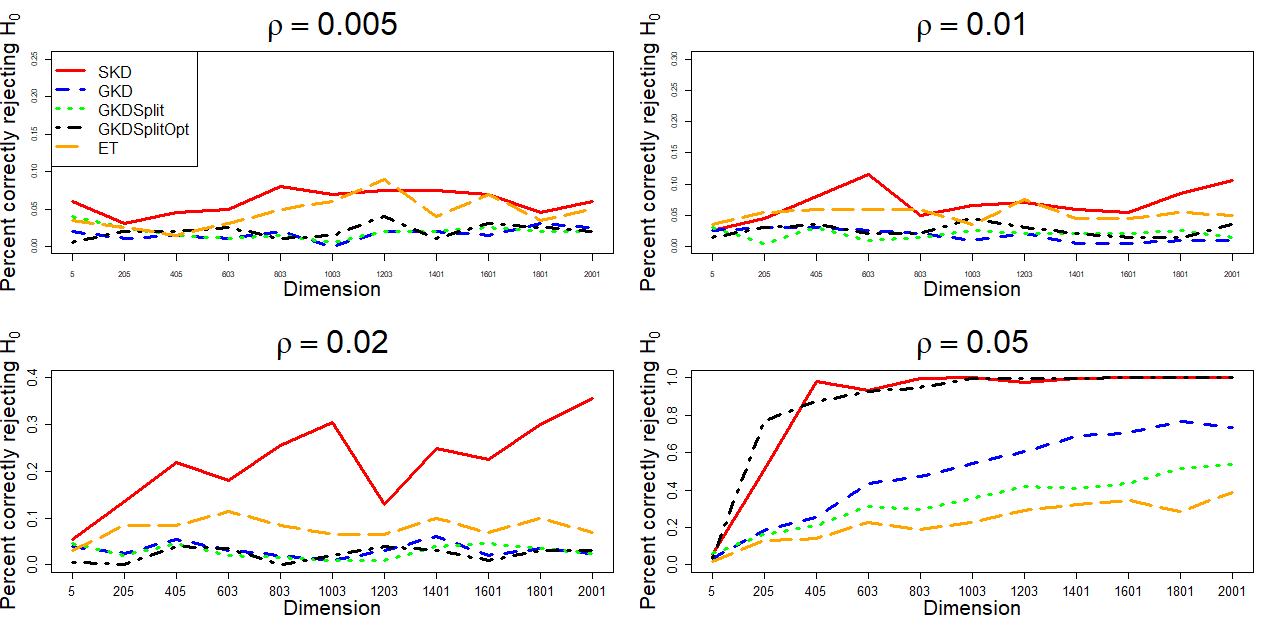}
					\caption{Performance of tests under Model 2.2 with $\alpha=0.05$. Four values of $\rho$ are shown: $0.005$, $0.01$, $0.02$ and $0.05$ (alternative hypothesis).}
					\label{fig:Model2.2}
				\end{figure}
				Results of Model 2.2 are shown in Fig. \ref{fig:Model2.2}. SKD seems to be particularly sensitive to dependence since correlations of $\rho=0.05$ quickly lead to a power of almost 1. SKD obtains the best results in this scenario.
			\subsubsection*{A real data example: Barcelona temperatures (1944-2019)}
				\begin{figure}[!h]
					\centering
					\includegraphics[width=0.9\linewidth]{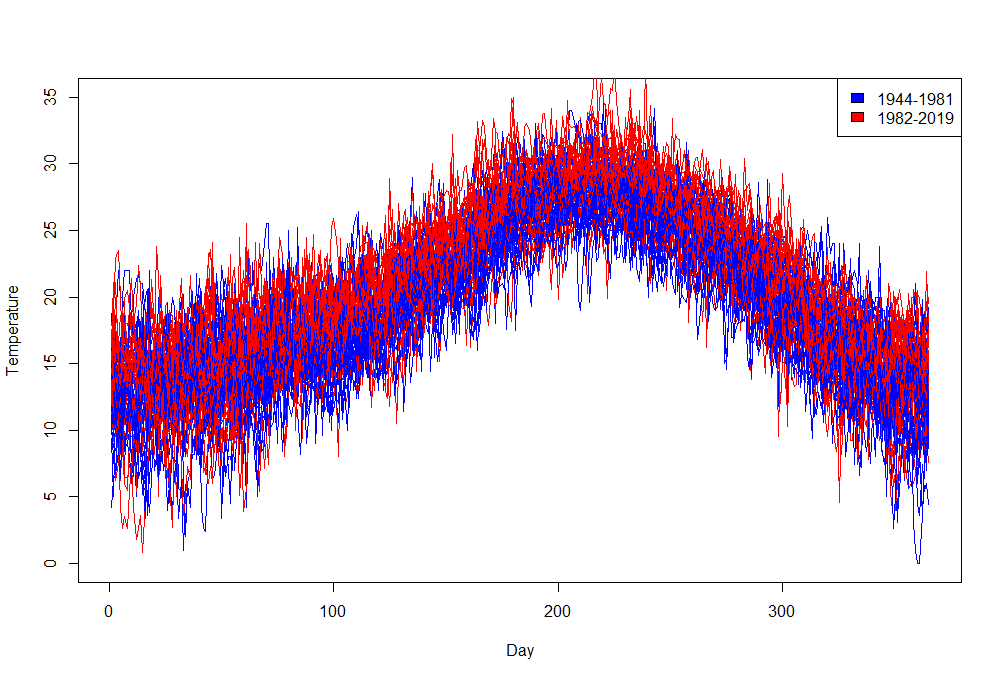}
					\caption{Maximum daily land surface temperature measured at El Prat Airport (Barcelona, Spain) between 1944 and 2019.}
					\label{fig:temperatures_BCN}
				\end{figure}
				We consider daily values of maximum temperatures registered at Barcelona airport (El Prat) from years 1944 to 2019. The data set consists of 76 vectors of dimension 365, each of which corresponds to a year in that time period. The daily observations have been treated as discretization points to include the problem within the framework of functional data, every year providing a function in the sample. Those observations corresponding to the 29th of February in leap years are omitted and missing observations are interpolated. These data are available at \url{https://www.ncei.noaa.gov}, the web page of the National Centers for Environmental Information.
				
				Our purpose is to test the null hypothesis that the sample of temperatures from 1944 to 1981 comes from the same (functional) distribution to that of the period 1982-2019. The rejection of this null hypothesis could be interpreted as a hint of possible warming in the area. Indeed, we observe that, in absence of any significant climate change, one would expect that both samples are made of independent trajectories from the same underlying process.

				All the considered tests give a nearly null $p$-value. This is hardly surprising, in view of Fig. \ref{fig:temperatures_BCN}, where the temperature curves are displayed (the blue curves correspond to the earlier period). While this is just a small experiment, presented here for illustration purposes, the results are consistent with those of many other deeper analysis published in recent years.
			\subsubsection*{Conclusions of the empirical study}
				In the light of the results, we can conclude that, globally, the supremum kernel distance test (SKD) performs similarly to the GKD tests in the homoscedastic case, though the ET test appears to be the winner in this situation. In the heteroscedastic case, SKD obtains almost the best power results. A more complete study (including the derivation of the asymptotic distribution for the case of a data-driven selection of $\lambda$, the use of Pearson curves and/or modified bootstrap schemes, \dots) might be worthwhile in the future. On the other hand, according to \cite{Fukumizu2013}, energy tests can be expressed in terms of kernel distances. This idea might deserve further attention as well, in order to incorporate these \enquote{equivalent} kernels to the SKD paradigm. In any case, it is clear that the present study does not allow us to conclude any obvious superiority (or inferiority) of none of the considered methods. In fact, the aim of our limited empirical study is to show that the SKD method can be implemented and it is competitive. This goal has been hopefully achieved. More definitive conclusions should be reached via subsequent empirical experiments and, especially, with the use of these tests by practitioners in the coming years. Software to run the SKD-based test will soon be available as an \verb!R!-package called \verb!SKD2!.
	\section{Final remarks}\label{sec:Concluding}
		The whole paper relies on the strategy of using kernel-based methods by combining, via supremum, different kernels. We consider the classical two-sample problem focusing on functional and high-dimensional data, where the demand for applicable methods seems more obvious. Despite the large amount of relevant literature on this topic, we believe that there is still room for improvement, as those provided here, in the line of obtaining more general results with a different technology of proofs. In particular, the use of differentiation techniques plus empirical processes methods allows us to address the asymptotic behaviour of the test statistics under the null and alternative hypothesis, including the case of unbalanced samples; see Theorems \ref{Theorem_asymptotic_H0} and \ref{Theorem_asymptotic_H1}. A key result in the proofs is Theorem \ref{Theorem:Donsker} that establishes, under some conditions covering infinite-dimensional settings, the universal Donsker property for a union of unit balls in RKHS. This theorem extends previous similar results for finite-dimensional situations and could be potentially useful in other statistical procedures within the RKHS framework. The approach established in this paper can also be potentially useful to analyze the asymptotic behavior of data-driven estimators of the kernel parameters. However, this interesting problem is beyond the scope of this paper. Theorems \ref{Theorem_asymptotic_H0} and \ref{Theorem_asymptotic_H1} are meaningful from a conceptual point of view, even if the limit distributions are not particularly simple. The mere existence of such limit (non-degenerate) distributions is a primary guarantee that kernel-based statistics can be used to derive procedures achieving, at least asymptotically, a prescribed significance level under the null and providing consistency under the alternative. We deal with the complicated structure of the limit distribution under the null by using permutation tests. This is also the case for other popular methodologies such as the energy test, \cite{Rizzo&Szekely2017}. Different approximation techniques are also conceivable, including truncation in the limit expression in Theorems \ref{Theorem_asymptotic_H0} plus estimation of the involved parameters.
While this paper was under review, other related and interesting contributions have appeared; see \cite{Albert2023} and \cite{Biggs2023}.

	\section{Technical details}\label{sec:Proofs}
		We need two auxiliary lemmata to prove Theorem \ref{Theorem:Donsker}. The first result corresponds to \cite[Theorem 1.1]{Marcus1985}.
		\begin{Lema}\label{Lemma.Markus}
			Let $H$ be a real and separable Hilbert space. Let us consider a linear and continuous operator $T:H\rightarrow\mathcal{C}_{b}(\mathcal{X})$, where $\mathcal{C}_{b}(\mathcal{X})$ is the space of real bounded continuous functions on $\mathcal{X}$ endowed with the supremum norm. If $B_{H}$ is the unit ball in $H$, then the class $B=T\left(B_{H}\right)$ is universal Donsker.
		\end{Lema}
		We also require Aronszajn's inclusion theorem; see \cite[Theorem I]{Aronszajn1950}.
		\begin{Lema}\label{Lemma.Aronszajn}
			Let $k_{1}$ and $k_{2}$ be two kernels on $\mathcal{X}$. Then, $\mathcal{H}_{k_{1}}\subset\mathcal{H}_{k_{2}}$ if and only if there exists a constant $c>0$ such that $c\,k_{2}-k_{1}$ is a positive definite kernel (i.e., $k_{1}\ll c\,k_{2}$). In such a case, we also have that $\|f\|_{\mathcal{H}_{k_{2}}}\leq\sqrt{c}\,\|f\|_{\mathcal{H}_{k_{1}}}$, for all $f\in\mathcal{H}_{k_{1}}$.
		\end{Lema}
		\begin{proof}[\textbf{\upshape Proof of Theorem \ref{Theorem:Donsker}}]
			The first part can be seen as a consequence of Lemma \ref{Lemma.Markus}. First, by \cite[Theorem 17]{Berlinet&Thomas-Agnan2011}, the functions in $\mathcal{H}_{k}$ are continuous. In particular, using \cite[Corollary 3]{Berlinet&Thomas-Agnan2011}, we conclude that $\mathcal{H}_{k}$ is a separable Hilbert space. On the other hand, for $x\in\mathcal{X}$, by the reproducing property \eqref{eq:reproducing} of $k$ (twice) and Cauchy–Schwarz inequality, we have that
			\begin{equation}\label{eq:inclusion_map}
				|f(x)-g(x)|=\left|\langle f-g,k(x,\cdot)\rangle_{\mathcal{H}_{k}}\right|\leq\| f-g\|_{\mathcal{H}_{k}}\,\|k(x,\cdot)\|_{\mathcal{H}_{k}}=\|f-g\|_{\mathcal{H}_{k}}\,\sqrt{k(x,x)}.
			\end{equation}
			Therefore, as $k$ is bounded on the diagonal, convergence in the RKHS norm entails uniform convergence. Further, from \eqref{eq:inclusion_map} we also see that the functions in $\mathcal{H}_{k}$ are bounded and hence $\mathcal{H}_{k}\subset\mathcal{C}_{b}(\mathcal{X})$.
			Now, we can apply Lemma \ref{Lemma.Markus} to $H=\mathcal{H}_{k}$ and $T=I$, the inclusion map given by $I(f)=f$. According to \eqref{eq:inclusion_map}, this linear transformation is continuous. As $B_{H}=\mathcal{F}_{k}$, by Lemma \ref{Lemma.Markus}, we thus conclude that $\mathcal{F}_{k}=T\left(\mathcal{F}_{k}\right)$ is universal Donsker.
			
			The second part is a by-product of the first one together with Aronszajn's inclusion theorem. According to Lemma \ref{Lemma.Aronszajn}, we have that $\|f\|_{\mathcal{H}_{k}}\leq\sqrt{c}\,\|f\|_{\mathcal{H}_{k,\lambda}}$, for all $f\in\mathcal{H}_{k,\lambda}$ and for all $\lambda\in\Lambda$. Therefore, $\mathcal{F}_{k,\Lambda}\subset \sqrt{c}\,\mathcal{F}_{k}$. Finally, from the first part of the theorem, the set $\sqrt{c}\,\mathcal{F}_{k}$ is universal Donsker as it is the unit ball of the RKHS generated by the kernel $ck$. Therefore, $\mathcal{F}_{k,\Lambda}$ is also universal Donsker (see \cite[Theorem 2.10.1]{van_der_Vaart&Wellner2023}) and the proof is complete.
		\end{proof}
		To prove Theorem \ref{Theorem_asymptotic_H0}, we need the following Karhunen-Loève-type result for the $\mathcal{F}_{k,\lambda}$-indexed Brownian bridge.
		\begin{Lema}\label{Lemma.KL}
			Under the assumptions of Theorem \ref{Theorem_asymptotic_H0}, we have that:
			\begin{description}
				\item[(a)]\label{item:Lemma.dual} For each $\lambda\in\Lambda$, the $\mathcal{F}_{k,\lambda}$-indexed Brownian bridge $\mathbb{G}_{\operatorname{P}}$ can be extended almost surely to a continuous and linear map on $\mathcal{H}_{k,\lambda}$. Therefore, $\mathbb{G}_{\operatorname{P}}$ can be seen as a random element of the dual space $\mathcal{H}^{\ast}_{k,\lambda}$. For simplicity we also denote this extension in $\mathcal{H}^{\ast}_{k,\lambda}$ as $\mathbb{G}_{\operatorname{P}}$.
				\item[(b)]\label{item:Lemma.KL} As an element of $\mathcal{H}^{\ast}_{k,\lambda}$, $\mathbb{G}_{\operatorname{P}}$ admits the following representation:
				\begin{equation}\label{KL-dual}
					\mathbb{G}_{\operatorname{P}}=_{\operatorname{a.s.}}\sum_{j\in\mathbb{N}}Z_{j,\lambda}\,\varphi_{j,\lambda},\quad\operatorname{in}\ \mathcal{H}_{k,\lambda}^{\ast}.
				\end{equation}
				In particular, we have that
				\begin{equation}\label{proof-representation}
					\left\|\mathbb{G}_{\operatorname{P}}\right\|_{\mathcal{H}_{k,\lambda}^{\ast}}^{2}=_{\operatorname{a.s.}}\sum_{j\in\mathbb{N}}Z_{j,\lambda}^{2}.
				\end{equation}
			\end{description}
		\end{Lema}
		\begin{proof}[\textbf{\upshape Proof:}]
			To show part \emph{(\ref{item:Lemma.dual})}, we note that, from Theorem \ref{Theorem:Donsker}, $\mathcal{F}_{k,\lambda}$ is a $\operatorname{P}$-Donsker class and hence $\operatorname{P}$-pre-Gaussian. Hence, by \cite[Theorem 3.7.28]{Gine&Nickl2016}, for almost all $\omega$, the function $f\mapsto\mathbb{G}_{\operatorname{P}}(\omega)f$ ($f\in\mathcal{F}_{k,\lambda}$) is prelinear and can be uniquely extended to a linear map on $\operatorname{span}\left(\mathcal{F}_{k,\lambda}\right)=\mathcal{H}_{k,\lambda}$. Moreover, this extension is bounded and uniformly $d_{\operatorname{P}}$-continuous in $\mathcal{H}_{k,\lambda}$, where $d_{\operatorname{P}}$ is the intrinsic $\operatorname{L}^{2}(\operatorname{P})$ metric of the process. Finally, we observe that, thanks to \eqref{eq:inclusion_map},
			\begin{equation}\label{eq:bound-L2-RKHS}
				d^{2}_{\operatorname{P}}(f,g)=\operatorname{E}_{\operatorname{P}}(f-g)^{2}\leq\|f-g\|_{\mathcal{H}_{k,\lambda}}^{2}\,\int_{\mathcal{X}}k(x,x)\,\operatorname{dP}(x).
			\end{equation}
			As by hypothesis $k$ is bounded on the diagonal, we have that uniformly $d_{\operatorname{P}}$-continuous functions on $\mathcal{H}_{k,\lambda}$ are also uniformly continuous functions with respect of the norm in $\mathcal{H}_{k,\lambda}$. In particular, $\mathbb{G}_{\operatorname{P}}$ is almost surely a continuous and linear functional on $\mathcal{H}_{k,\lambda}$, and thus an element of $\mathcal{H}^{\ast}_{k,\lambda}$. This finishes the proof of part \emph{(\ref{item:Lemma.dual})}.

			To prove part \emph{(\ref{item:Lemma.KL})} we first note that, by  \emph{(\ref{item:Lemma.dual})}, $\mathbb{G}_{\operatorname{P}}$ is a Gaussian process in the Hilbert space $\mathcal{H}_{k,\lambda}^{\ast}$. The covariance operator $\mathcal{K}_{\mathbb{G}_{\operatorname{P}}}$ of $\mathbb{G}_{\operatorname{P}}$ is self-adjoint and compact. By the Fernique's theorem \cite[p. 74]{Bogachev1998}, $\mathbb{G}_{\operatorname{P}}$ is Bochner square-integrable, $\mathcal{K}_{\mathbb{G}_{\operatorname{P}}}$ is a trace-class operator and
			\begin{equation}\label{formula: trace}
				\operatorname{trace}\left(\mathcal{K}_{\mathbb{G}_{\operatorname{P}}}\right)=\int_{\mathcal{H}_{k,\lambda}^{\ast}}\|z\|_{\mathcal{H}_{k,\lambda}^{\ast}}^{2}\,\operatorname{d}\nu_{\mathbb{G}_{\operatorname{P}}}(z)=\mathbb{E}\left(\left\|\mathbb{G}_{\operatorname{P}}\right\|_{\mathcal{H}_{k,\lambda}^{\ast}}^{2}\right),
			\end{equation}
			where $\nu_{\mathbb{G}_{\operatorname{P}}}$ is the measure induced by the process $\mathbb{G}_{\operatorname{P}}$ in $\mathcal{H}_{k,\lambda}^{\ast}$. The proof of \eqref{formula: trace} can be found in \cite[p. 48]{Bogachev1998}.

			Now, by the spectral theorem, there exists $\left\{\left(\beta_{j,\lambda},\varphi_{j,\lambda}\right)\right\}_{j\in\mathbb{N}}\in\left([0,\infty)\times\mathcal{H}_{k,\lambda}^{\ast}\right)^{\mathbb{N}}$ such that $\beta_{1,\lambda}\geq\beta_{2,\lambda}\geq\cdots$; $\mathcal{K}_{\mathbb{G}_{P}}\,\varphi_{j,\lambda}=\beta_{j,\lambda}\,\varphi_{j,\lambda}$, for $j\in\mathbb{N}$; and $\left\langle\varphi_{j_{1},\lambda},\varphi_{j_{2},\lambda}\right\rangle_{\mathcal{H}_{k,\lambda}^{\ast}}=\delta_{j_{1}\,j_{2}}$, for $j_{1},j_{2}\in\mathbb{N}$ with $\delta_{i\,j}$ the Kronecker's delta. As $\mathcal{K}_{\mathbb{G}_{\operatorname{P}}}$ is trace-class, we also have that $\operatorname{trace}\left(\mathcal{K}_{\mathbb{G}_{P}}\right)=\sum_{j\in\mathbb{N}}\beta_{j,\lambda}.$ Additionally,
			\begin{equation}\label{Proof-normal}
				\operatorname{E}\left(\left\langle\mathbb{G}_{\operatorname{P}},\varphi_{j,\lambda}\right\rangle_{\mathcal{H}_{k,\lambda}^{\ast}}\right)=0\quad\operatorname{and}\quad\operatorname{E}\left(\left\langle\mathbb{G}_{\operatorname{P}},\varphi_{j,\lambda}\right\rangle_{\mathcal{H}_{k,\lambda}^{\ast}}^{2}\right)=\left\langle\mathcal{K}_{\mathbb{G}_{\operatorname{P}}}\left(\varphi_{j,\lambda}\right),\varphi_{j,\lambda}\right\rangle_{\mathcal{H}_{k,\lambda}^{\ast}}=\beta_{j,\lambda}.
			\end{equation}
			From \eqref{Proof-normal}, we have that $Z_{j,\lambda}=\left\langle\mathbb{G}_{\operatorname{P}},\varphi_{j,\lambda}\right\rangle_{\mathcal{H}_{k,\lambda}^{\ast}}\sim\mathcal{N}\left(0,\beta_{j,\lambda}\right)$ ($j\in\mathbb{N}$) are jointly Gaussian and independent.
			
			To finish this proof of \eqref{KL-dual}, by \cite[Theorem 6.1]{Ledoux&Talagrand2013}, it is enough to show absolute mean convergence, which is a necessary and sufficient condition. First, by orthogonality, we observe that for every $J\subset\mathbb{N}$ finite, we have that
			\begin{equation*}
				0\leq\left\|\mathbb{G}_{\operatorname{P}}-\sum_{j\in J}Z_{j,\lambda}\,\varphi_{j,\lambda}\right\|_{\mathcal{H}_{k,\lambda}^{\ast}}^{2}=\left\|\mathbb{G}_{\operatorname{P}}\right\|_{\mathcal{H}_{k,\lambda}^{\ast}}^{2}-\sum_{j\in J}Z_{j,\lambda}^{2}.
			\end{equation*}
			Then by \eqref{formula: trace},
			\begin{equation}\label{formula: series_remainder}
				\operatorname{E}\left(\left|\left\|\mathbb{G}_{\operatorname{P}}\right\|_{\mathcal{H}_{k,\lambda}^{\ast}}^{2}-\sum_{j\in J}Z_{j,\lambda}^{2}\right|\right)=\operatorname{trace}\left(\mathcal{K}_{\mathbb{G}_{\operatorname{P}}}\right)-\sum_{j\in J}\beta_{j,\lambda},
			\end{equation}
			which is the remainder of a convergent series. Hence, \eqref{KL-dual} holds. As \eqref{proof-representation} follows from \eqref{KL-dual}, the proof is complete.
		\end{proof}
		The proof of part \emph{(\ref{item:Lemma.dual})} in Lemma \ref{Lemma.KL} essentially follows from Theorem \ref{Theorem:Donsker}. However, part \emph{(\ref{item:Lemma.KL})}, where the series representation is obtained, must be discussed. Equation \eqref{KL-dual} shows the convergence of a series of functional random variables. This result looks like a standard Karhunen-Loève theorem, but some remarks should be done. The convergence of this series is on the dual space $\mathcal{H}_{k,\lambda}^{\ast}$, while Karhunen-Loève decomposition is stated classically on $\operatorname{L}^{2}$-type spaces. In fact, our decomposition in \eqref{KL-dual} can be seen as a particular case of the results in \cite{Bay&Croix2017}. In \cite[Theorem 2.6.10]{Gine&Nickl2016} a similar decomposition is shown where the coordinates are deterministic while the basis is random, which is not useful for our purposes.
		\begin{proof}[\textbf{\upshape Proof of Theorem \ref{Theorem_asymptotic_H0}}]
			From Theorem \ref{Theorem:Donsker}, the class $\mathcal{F}_{k,\Lambda}$ is Donsker and hence we have that
			\begin{equation}\label{eq:Th2-1}
				\mathbb{G}_{n,m}=\sqrt{\frac{n\,m}{n+m}}\,\left(\mathbb{P}_{n}-\mathbb{Q}_{m}\right)\rightsquigarrow\mathbb{G}_{\operatorname{P}},\quad\operatorname{in}\ \ell^{\infty}\left(\mathcal{F}_{k,\Lambda}\right).
			\end{equation}
			Note that $d_{k,\Lambda}$ is the metric induced by the supremum norm in $\ell^{\infty}\left(\mathcal{F}_{k,\Lambda}\right)$, hence $d_{k,\Lambda}$ is a continuous functional. From \eqref{eq:Th2-1} and by the continuous mapping theorem (see, for instance \cite[Theorem 1.9.5]{van_der_Vaart&Wellner2023}), we obtain that
			\begin{equation}\label{eq:Th2-2}
				\sqrt{\frac{n\,m}{n+m}}\,d_{k,\Lambda}\left(\mathbb{P}_{n},\mathbb{Q}_{m}\right)\rightsquigarrow\underset{\mathcal{F}_{k,\Lambda}}{\operatorname{sup}}\left(\mathbb{G}_{\operatorname{P}}\right).
			\end{equation}
			From Lemma \ref{Lemma.KL}, the limit in \eqref{eq:Th2-2} can be rewritten as
			\begin{equation}\label{eq:Th2-3}
				\underset{\mathcal{F}_{k,\Lambda}}{\operatorname{sup}}\,\left(\mathbb{G}_{\operatorname{P}}\right)=\underset{\lambda\in\Lambda}{\operatorname{sup}}\,\left(\underset{\mathcal{F}_{k,\lambda}}{\operatorname{sup}}\,\left(\mathbb{G}_{\operatorname{P}}\right)\right)=\underset{\lambda\in\Lambda}{\operatorname{sup}}\,\left(\left\|\mathbb{G}_{\operatorname{P}}\right\|_{\mathcal{H}_{k,\lambda}^{\ast}}\right).
			\end{equation}
			Finally, from \eqref{proof-representation} and \eqref{eq:Th2-3} we obtain the representation of the limit as in \eqref{eq:theorem-H0} and the proof of the theorem is complete.
		\end{proof}
		The next goal is to prove Corollary \ref{Corollary-H1} as preparation for the proof of Theorem \ref{Theorem_asymptotic_H1}. We need a differentiability result for the supremum similar to those obtained in \cite{Carcamo2020}. Given a kernel $k$, we consider the mapping
		\begin{equation}\label{eq:sigma}
			\sigma_{k}(g)=\underset{f\in \mathcal{F}_{k}}{\operatorname{sup}}\,(g(f)),\quad\ \operatorname{for}\ g\in\ell^{\infty}\left(\mathcal{F}_{k}\right),
		\end{equation}
		where $\mathcal{F}_{k}$ is the unit ball as in (\ref{unit-ball}). Observe that, by \eqref{definition.kernel.distance} and \eqref{kernel.integral-metric}, if $\operatorname{P},\operatorname{Q}\in\mathcal{M}_{\operatorname{p}}(\mathcal{X})$ such that their mean embeddings $\mu_{\operatorname{P}}$ and $\mu_{\operatorname{Q}}$ exist, we have that
		\begin{equation}\label{eq:sigma-kd}
			\sigma_{k}(\operatorname{P}-\operatorname{Q})=d_{k}(\operatorname{P},\operatorname{Q})=\left\|\mu_{\operatorname{P}}-\mu_{\operatorname{Q}}\right\|_{\mathcal{H}_{k}}.
		\end{equation}
		The proof of Corollary \ref{Corollary-H1} relies on Theorem \ref{Theorem:Donsker} together with the differentiability properties of the mapping $\sigma_{k}$ in \eqref{eq:sigma}. Same ideas are used below in the proof of Theorem \ref{Theorem_asymptotic_H1} using the mapping
		\begin{equation}\label{eq:sigma-Lambda}
			\sigma_{k,\Lambda}(g)=\underset{f\in \mathcal{F}_{k,\Lambda}}{\operatorname{sup}}\,(g(f)),\quad  \operatorname{for}\ g\in\ell^{\infty}\left(\mathcal{F}_{k,\Lambda}\right),
		\end{equation}
		where $\mathcal{F}_{k,\Lambda}$ is the union of balls in \eqref{SKD-integral}. These differentiability results might have independent interest as it can be applied in other contexts by means of the (extended) functional Delta method; see the examples in \cite{Carcamo2020}.
		
		The next corollary shows that $\sigma_{k}$ in \eqref{eq:sigma} is fully Hadamard differentiable under some assumptions. For the precise definitions we refer to \cite{Carcamo2020} and the references therein.
		\begin{Lema}\label{Lemma-Differentiability-1}
			Let us consider $\operatorname{P},\operatorname{Q}\in\mathcal{M}_{\operatorname{p}}(\mathcal{X})$ such that their mean embeddings $\mu_{\operatorname{P}}$ and $\mu_{\operatorname{Q}}$ exist and $\mu_{\operatorname{P}}\neq\mu_{\operatorname{Q}}$. We have that the mapping $\sigma_{k}$ in \eqref{eq:sigma} is (fully) Hadamard differentiable at $\operatorname{P}-\operatorname{Q}$ tangentially to $\mathcal{C}\left(\mathcal{F}_{k},d_{\mathcal{H}_{k}}\right)\equiv$ the subset of $\ell^{\infty}\left(\mathcal{F}_{k}\right)$ constituted by continuous functions with respect to
			the RKHS norm. In such a case, the derivative of $\sigma_{k}$ at the point $\operatorname{P}-\operatorname{Q}$ is given by
			\begin{equation}\label{sigma.prime.RKHS}
				\sigma_{k;\operatorname{P}-\operatorname{Q}}^{\prime}(g)= g\left(h^{+}\right),\quad\operatorname{for}\ g\in\mathcal{C}\left(\mathcal{F}_{k},d_{\mathcal{H}_{k}}\right),
			\end{equation}
			where $h^{+}\in\mathcal{F}_{k}$ is defined in \eqref{h+funtion}.
		\end{Lema}
		\begin{proof}[\textbf{\upshape Proof:}]
			From \cite[Theorem 2.1]{Carcamo2020}, we have that $\sigma_{k}$ is Hadamard directionally differentiable and
			\begin{equation}\label{eq:derivative-1}
				\sigma_{k;\operatorname{P}-\operatorname{Q}}^{\prime}(g)=\underset{\varepsilon\searrow0}{\operatorname{lim}}\,\underset{f\in A_{\varepsilon}(\operatorname{P}-\operatorname{Q})}{\operatorname{sup}}\,(g(f)),\quad g\in\ell^{\infty}\left(\mathcal{F}_{k}\right),
			\end{equation}
			where
			$A_{\varepsilon}(\operatorname{P}-\operatorname{Q})=\left\{h\in\mathcal{F}_{k}:(\operatorname{P}-\operatorname{Q})\,(h)\geq d_{k}(\operatorname{P},\operatorname{Q})-\varepsilon\right\}.$ We first check that if $h_{\varepsilon}\in A_{\varepsilon}(\operatorname{P}-\operatorname{Q})$, then $h_{\varepsilon}\to h^{+}$ in $\mathcal{H}_{k}$ as $\varepsilon\to0$, with $h^{+}$ in \eqref{h+funtion}. To see this, we first note that
			\begin{equation}\label{h.epsion}
				\left\|h_{\varepsilon}-h^{+}\right\|_{\mathcal{H}_{k}}^{2}=1+\left\|h_{\varepsilon}\right\|_{\mathcal{H}_{k}}^{2}-\frac{2}{\left\|\mu_{\operatorname{P}}-\mu_{\operatorname{Q}}\right\|_{\mathcal{H}_k}}\,\inp{h_{\varepsilon}}{\mu_{\operatorname{P}}-\mu_{\operatorname{Q}}}{\mathcal{H}_{k}}.
			\end{equation}
			As $h_{\varepsilon}\in A_{\varepsilon}(\operatorname{P}-\operatorname{Q})$, from \eqref{eq:mean-embedding-property} and (\ref{eq:sigma-kd}), we obtain that
			\begin{equation}\label{eq:dif.bound}
				\operatorname{P}\left(h_{\varepsilon}\right)-\operatorname{Q}\left(h_{\varepsilon}\right)=\inp{h_{\varepsilon}}{\mu_{\operatorname{P}}-\mu_{\operatorname{Q}}}{\mathcal{H}_{k}}\geq\left\|\mu_{\operatorname{P}}-\mu_{\operatorname{Q}}\right\|_{\mathcal{H}_{k}}-\varepsilon.
			\end{equation}
			Finally, from (\ref{h.epsion}), \eqref{eq:dif.bound}, and as $h_{\varepsilon}\in\mathcal{F}_{k}$, we have that
			\begin{equation}\label{eq:dif.bound-2}
				\left\|h_{\varepsilon}-h^{+}\right\|_{\mathcal{H}_{k}}^{2}\leq\frac{2\,\varepsilon}{\left\|\mu_{\operatorname{P}}-\mu_{\operatorname{Q}}\right\|_{\mathcal{H}_{k}}},
			\end{equation}
			and hence $h_{\varepsilon}\to h^{+}$ in $\mathcal{H}_{k}$ (as $\varepsilon\to0$).
			
			Now, we check that $\sigma_{k;\operatorname{P}-\operatorname{Q}}^{\prime}(g)=g\left(h^{+}\right)$, for $g\in\mathcal{C}\left(\mathcal{F}_{k},d_{\mathcal{H}_{k}}\right)$. We firstly observe that $h^{+}\in A_{\varepsilon}(\operatorname{P}-\operatorname{Q})$, for all $\varepsilon>0$. Hence, from equation (\ref{eq:derivative-1}), we have that $g\left(h^{+}\right)\leq\sigma_{k;\operatorname{P}-\operatorname{Q}}^{\prime}(g)$. On the other hand, we can extract a maximizing sequence $h_{m}\in A_{\rfrac{1}{m}}(\operatorname{P}-\operatorname{Q})$ ($m\in\mathbb{N}$) satisfying that $\underset{A_{\rfrac{1}{m}}(\operatorname{P}-\operatorname{Q})}{\operatorname{sup}}\,(g)\leq g\left(h_{m}\right)+\frac{1}{m}$. As $g$ is continuous and $h_{m}\to h^{+}$ as $m\to\infty$ in $\mathcal{H}_{k}$, we obtain that $\sigma_{k;\operatorname{P}-\operatorname{Q}}^{\prime}(g)=\underset{m\to\infty}{\operatorname{lim}}\,\underset{A_{\rfrac{1}{m}}(\operatorname{P}-\operatorname{Q})}{\operatorname{sup}}\,(g)\leq\underset{m\to\infty}{\operatorname{lim}}\,g\left(h_{m}\right)=g\left(h^{+}\right).$ Therefore, we obtain that $\sigma{k;\operatorname{P}-\operatorname{Q}}^{\prime}(g)=g\left(h^{+}\right)$, which is a linear mapping, so $\sigma_{k}$ is fully differentiable and the proof is complete.
		\end{proof}
		\begin{proof}[Proof of Corollary \ref{Corollary-H1}]
			From Theorem \ref{Theorem:Donsker}, we have that
			\begin{equation}\label{eq:limit-G}
				\mathbb{G}_{n,m}=\sqrt{\frac{n\,m}{n+m}}\,\left(\mathbb{P}_{n}-\mathbb{Q}_{m}-(\operatorname{P}-\operatorname{Q})\right)\rightsquigarrow\mathbb{G}=\sqrt{1-\theta}\,\mathbb{G}_{\operatorname{P}}-\sqrt{\theta}\,\mathbb{G}_{\operatorname{Q}},\quad\operatorname{in}\ \ell^{\infty}\left(\mathcal{F}_{k}\right).
			\end{equation}
			From \eqref{eq:sigma-kd}, the statistic in the right-hand side of equation \eqref{eq:corolary1} is precisely
			\begin{equation}
				\sqrt{\frac{n\,m}{n+m}}\,\left(\sigma_{k}\left(\mathbb{P}_{n}-\mathbb{Q}_{m}\right)-\sigma_{k}(\operatorname{P}-\operatorname{Q})\right),
			\end{equation}
			where $\sigma_{k}$ is defined in \eqref{eq:sigma}.
			
			Using the same ideas as in the proof of \cite[Theorem 6.1]{Carcamo2020}, it can be checked that the paths of $\mathbb{G}$ in \eqref{eq:limit-G} are a.s. in $\mathcal{C}_{u}\left(\mathcal{F}_{k},\rho\right)$ (uniformly continuous), where
			\begin{equation}\label{eq:metric-rho}
				\rho=\operatorname{max}\left(d_{\operatorname{L}^{2}(\operatorname{P})},d_{\operatorname{L}^{2}(\operatorname{Q})}\right),
			\end{equation}
			is the natural $\operatorname{L}^{2}$-metric of $\mathbb{G}$. From \eqref{eq:bound-L2-RKHS}, it can be readily checked that $\mathcal{C}_{u}\left(\mathcal{F}_{k},\rho\right)\subset\mathcal{C}_{u}\left(\mathcal{F}_{k},d_{\mathcal{H}_{k}}\right)$ and hence $\mathbb{G}\in\mathcal{C}\left(\mathcal{F}_{k},d_{\mathcal{H}_{k}}\right)$ a.s. To finish the proof it is enough to apply Lemma \ref{Lemma-Differentiability-1} together with the functional Delta method \cite[Section 3.10]{van_der_Vaart&Wellner2023}.
		\end{proof}
		To prove Theorem \ref{Theorem_asymptotic_H1} we need the following key lemma.
		\begin{Lema}\label{Lemma-Differentiability-2}
			Let us assume that the family of kernels $\left\{k_{\lambda}:\lambda\in\Lambda\right\}$ satisfies \ref{itm:Dom}, \ref{itm:Ide} and \ref{itm:Par}. If $\operatorname{P},\operatorname{Q}\in\mathcal{M}_{\operatorname{p}}(\mathcal{X})$ such that $\operatorname{P}\neq\operatorname{Q}$, then the mapping $\sigma_{k,\Lambda}$ in \eqref{eq:sigma-Lambda} is Hadamard directionally differentiable at $\operatorname{P}-\operatorname{Q}$ tangentially to $\mathcal{C}\left(\mathcal{F}_{k,\Lambda},\rho\right)\equiv$ the subset of $\ell^{\infty}\left(\mathcal{F}_{k,\Lambda}\right)$ constituted by continuous functionals with respect to the distance $\rho$ in \eqref{eq:metric-rho}. In such a case, the (directional) derivative of $\sigma_{k,\Lambda}$ at the point $\operatorname{P}-\operatorname{Q}$ is given by
			\begin{equation}\label{sigma.prime.RKHS-Lambda}
				\sigma_{k,\Lambda;\operatorname{P}-\operatorname{Q}}^{\prime}(g)=\underset{\lambda\in\Lambda_{0}}{\operatorname{sup}}\,\left(g\left(h^{+,\lambda}\right)\right)=\underset{L}{\operatorname{sup}}\,(g),\quad g\in\mathcal{C}\left(\mathcal{F}_{k,\Lambda},\rho\right),
			\end{equation}
			where the functions $h^{+,\lambda}$ are defined in \eqref{Witness} and the sets $\Lambda_{0}$ and $L$ in \eqref{Witness-sets}.
		\end{Lema}
		\begin{proof}[\textbf{\upshape Proof:}]
			Let us fix $g\in\mathcal{C}\left(\mathcal{F}_{k,\Lambda},\rho\right)$. Again, from \cite[Theorem 2.1]{Carcamo2020}, we have that $\sigma_{k,\Lambda}$ is Hadamard directionally differentiable and
			\begin{equation}
				\sigma_{k,\Lambda;\operatorname{P}-\operatorname{Q}}^{\prime}(g)=\underset{\varepsilon\searrow0}{\operatorname{lim}}\,\underset{A_{\varepsilon,\Lambda}(\operatorname{P}-\operatorname{Q})}{\operatorname{sup}}\,(g),
			\end{equation}
			where
			$A_{\varepsilon,\Lambda}(\operatorname{P}-\operatorname{Q})=\left\{h\in\mathcal{F}_{k,\Lambda}:(\operatorname{P}-\operatorname{Q})\,(h)\geq d_{k,\Lambda}(\operatorname{P},\operatorname{Q})-\varepsilon\right\}.$ For every $\varepsilon>0$, it is clear that $L\subseteq A_{\varepsilon,\Lambda}(\operatorname{P}-\operatorname{Q})$, where $L$ is defined in \eqref{Witness-sets}. Hence, we have that
			\begin{equation}\label{eq:dif-sup-bound-1}
				\underset{\lambda\in\Lambda_{0}}{\operatorname{sup}}\left(g\left(h^{+,\lambda}\right)\right)\leq\sigma_{k,\Lambda;\operatorname{P}-\operatorname{Q}}^{\prime}(g).
			\end{equation}
			Conversely, we consider a maximizing sequence $\left(h_{m}\right)_{m\in\mathbb{N}}$ satisfying that $h_{m}\in A_{\rfrac{1}{m},\Lambda}(\operatorname{P}-\operatorname{Q})$ and
			\begin{equation}\label{eq:Aset2}
				\underset{A_{\rfrac{1}{m},\Lambda}(\operatorname{P}-\operatorname{Q})}{\operatorname{sup}}\,(g)\leq g\left(h_{m}\right)+\frac{1}{m}.
			\end{equation}
			Each $h_{m}\in\mathcal{F}_{k,\lambda_{m}}$ ($m\in\mathbb{N}$), for some $\lambda_{m}\in\Lambda$. We consider the sequence $\left(h^{+,\lambda_{m}}\right)_{m\in\mathbb{N}}$. Using \ref{itm:Par}, by restricting, if needed, to a subsequence we can assume that $\lambda_{m}\to\lambda^{\star}\in\Lambda$. Next we prove the following facts:
			\begin{description}
				\labitem{(i)}{item:supremum.i} $\lambda^{\star}\in\Lambda_{0}$, where $\Lambda_{0}$ is in \eqref{Witness-sets}, and
				\begin{equation}\label{eq:sd-1}
					\left\|\mu_{\operatorname{P}}^{\lambda_{m}}-\mu_{\operatorname{Q}}^{\lambda_{m}}\right\|_{\mathcal{H}_{k,\lambda_{m}}}\to\left\|\mu_{\operatorname{P}}^{\lambda^{\star}}-\mu_{\operatorname{Q}}^{\lambda^{\star}}\right\|_{\mathcal{H}_{k,\lambda^{\star}}}=d_{k,\Lambda}(\operatorname{P},\operatorname{Q})>0.
				\end{equation}
				\labitem{(ii)}{item:supremum.ii}  $\rho\left(h_{m},h^{+,\lambda_{m}}\right)\to0$, as $m\to\infty$.
				\labitem{(iii)}{item:supremum.iii} $\rho\left(h_{m},h^{+,\lambda^{\star}}\right)\to0$, as $m\to\infty$.
			\end{description}
			First, \eqref{eq:sd-1} is obtained by using the representation of the kernel distance as a double integral in \eqref{definition.kernel.distance}, together with \ref{itm:Dom}, \ref{itm:Par} and the dominated convergence theorem (DCT). Further, as $h_{m}\in\mathcal{F}_{k,\lambda_{m}}\cap A_{\rfrac{1}{m},\Lambda}(\operatorname{P}-\operatorname{Q})$, we obtain that
			$$\left\|\mu_{\operatorname{P}}^{\lambda_{m}}-\mu_{\operatorname{Q}}^{\lambda_{m}}\right\|_{\mathcal{H}_{k,\lambda_{m}}}\geq\operatorname{P}\left(h_{m}\right)-\operatorname{Q}\left(h_{m}\right)\geq d_{k,\Lambda}(\operatorname{P},\operatorname{Q})-\frac{1}{m}.$$
			Hence, from \eqref{eq:sd-1} and by taking $m\to\infty$ we obtain that $\lambda^{\star}\in\Lambda_{0}$. The fact that $d_{k,\Lambda}(\operatorname{P},\operatorname{Q})>0$ follows from \ref{itm:Ide} and the proof of \ref{item:supremum.i} is complete.

			To show \ref{item:supremum.ii}, using the same ideas as in the proof of equation \eqref{eq:dif.bound-2} and \ref{item:supremum.i}, we obtain that
			\begin{equation*}
				\left\|h_{m}-h^{+,\lambda_{m}}\right\|^{2}_{\mathcal{H}_{k,\lambda_{m}}}\leq2-\frac{2\,d_{k,\Lambda}(\operatorname{P},\operatorname{Q})-\frac{1}{m}}{\left\|\mu_{\operatorname{P}}^{\lambda_{m}}-\mu_{\operatorname{Q}}^{\lambda_{m}}\right\|_{\mathcal{H}_{k,\lambda_{m}}}}\to0,\quad\operatorname{as}\ m\to\infty.
			\end{equation*}
			Now, from \eqref{eq:bound-L2-RKHS} and \ref{itm:Dom}, we have that for $\operatorname{S}\in\{\operatorname{P},\operatorname{Q}\}$
			$$d_{\operatorname{S}}^{2}\left(h_{m},h^{+,\lambda_{m}}\right)^{2}\leq\left\|h_{m}-h^{+,\lambda_{m}}\right\|_{\mathcal{H}_{k,\lambda_{m}}}^{2}\,c\,\int_{\mathcal{X}}k(x,x)\,\operatorname{dS}(x)\to0,\quad\operatorname{as}\ m\to\infty.$$
			Therefore, $\rho\left(h_{m},h^{+,\lambda_{m}}\right)\to0$, and \ref{item:supremum.ii} holds.
			
			To check \ref{item:supremum.iii},  by \ref{item:supremum.ii}, it is enough to see that $\rho\left(h^{+,\lambda_{m}},h^{+,\lambda^{\star}}\right)\to0$, as $m\to\infty$. By (\ref{eq:sd-1}) and repeatedly applying DCT (thanks to \ref{itm:Dom}), it can be checked that
			$$h^{+,\lambda_{m}}(x)\to h^{+,\lambda^{\star}}(x),\quad\operatorname{as}\ m\to\infty\ \operatorname{and\ for\ all}\ x\in\mathcal{X}.$$
			Furthermore, for $m$ large enough, we have that for $\operatorname{S}\in\{\operatorname{P},\operatorname{Q}\}$
			$$\left|h^{+,\lambda_{m}}(x)\right|\leq\frac{2\,c\,\left(\left|\mu_{\operatorname{P}}(x)\right|+\left|\mu_{\operatorname{Q}}(x)\right|\right)}{d_{k,\Lambda}(\operatorname{P},\operatorname{Q})}\in\operatorname{L}^{2}(\operatorname{S}),$$
			where $\mu_{\operatorname{P}}$ and $\mu_{\operatorname{Q}}$ are the mean embeddings corresponding to the dominating kernel $k$ in \ref{itm:Dom}. Hence, we can apply one more time DCT to obtain that $d_{\operatorname{S}}\left(h^{+,\lambda_{m}},h^{+,\lambda^{\star}}\right)\to0$. This implies that $\rho\left(h^{+,\lambda_{m}},h^{+,\lambda^{\star}}\right)\to0$, as $m\to\infty$.

			To finish, we use \eqref{eq:Aset2}, \ref{item:supremum.i}, \ref{item:supremum.iii}, as well as the continuity of the functional $g$ (with respect to the metric $\rho$) to obtain that
			\begin{equation*}
				\sigma_{k,\Lambda;\operatorname{P}-\operatorname{Q}}^{\prime}(g)=\underset{m\to\infty}{\operatorname{lim}}\,\underset{A_{\rfrac{1}{m},\Lambda}(\operatorname{P}-\operatorname{Q})}{\operatorname{sup}}\,(g)\leq\underset{m\rightarrow\infty}{\operatorname{lim}}\,g\left(h_{m}\right)=g\left(h^{+,\lambda^{\star}}\right)\leq\underset{\lambda\in\Lambda_{0}}{\operatorname{sup}}\left(g\left(h^{+,\lambda}\right)\right)=\underset{L}{\operatorname{sup}}\,(g).
			\end{equation*}
			The conclusion of this lemma follows from \eqref{eq:dif-sup-bound-1} and the previous inequalities.
		\end{proof}
		\begin{proof}[Proof of Theorem \ref{Theorem_asymptotic_H1}]
			The proof of this theorem is analogous to that of Corollary \ref{Corollary-H1} using Lemma \ref{Lemma-Differentiability-2} instead of Lemma \ref{Lemma-Differentiability-1}. Details are ommited.
		\end{proof}
	\section*{Acknowledgments}
This work has been partially supported by Grant PID2019-109387GB-I00 from the Spanish Ministry of Science and Innovation. A. Cuevas acknowledges as well financial support from Grant CEX2019-000904-S funded
by MCIN/AEI/10.13039/501100011033. L.A. Rodr\'{i}guez acknowledges funding from the German Science Foundation DFG Research Unit 5381: ``Mathematical Statistics in the Information Age–Statistical Efficiency and Computational Tractability''. The authors would like to thank two anonymous reviewers for their comments and references on the first version of this paper. We are also grateful to George Wynne (Imperial College, London) for useful comments regarding the dominance hypothesis for families of kernels in the infinite-dimensional case. We also want to thank Bojan Mihaljevic (Universidad Polit\'ecnica de Madrid) for his valuable assistance and advice in the development of the empirical section. Finally, kind help from Francisco de la Hoz (University of the Basque Country) with the final computations is also gratefully acknowledged. Some of the computational work has been done in the Center for Scientific Computation (UAM, project \textit{estfunc}).

		%

\begin{thebibliography}{99}
		\bibitem[{Albert et al.(2023)}]{Albert2023}
		\bibinfo{author}{M.~Albert}, \bibinfo{author}{B.~Guedj}, \bibinfo{author}{A.~Gretton}, \bibinfo{author}{I.~Kim}, \bibinfo{author}{W.~Xu}, \bibinfo{author}{A.~Schrab}, \bibinfo{title}{MMD aggregated two-sample test}, \bibinfo{journal}{Journal of Machine Learning Research}, \bibinfo{volume}{24} (\bibinfo{year}{2023}) \bibinfo{pages}{1--81}.
		
		\bibitem[{Aronszajn(1950)}]{Aronszajn1950}
		\bibinfo{author}{N.~Aronszajn}, \bibinfo{title}{Theory of reproducing kernels}, \bibinfo{journal}{Transactions of the American mathematical society}, \bibinfo{volume}{68} (\bibinfo{year}{1950}) \bibinfo{pages}{337--404}.
		
		\bibitem[{Balakrishnan et al.(2012)}]{Balakrishnan2012}
		\bibinfo{author}{S.~Balakrishnan}, \bibinfo{author}{K.~Fukumizu}, \bibinfo{author}{A.~Gretton}, \bibinfo{author}{M.~Pontil}, \bibinfo{author}{D.~Sejdinovic}, \bibinfo{author}{B.~K.~Sriperumbudur}, \bibinfo{author}{H.~Strathmann}, \bibinfo{title}{Optimal kernel choice for large-scale two-sample tests}, \bibinfo{journal}{Advances in neural information processing systems}, \bibinfo{volume}{25} (\bibinfo{year}{2012}).
		
		\bibitem[{del Barrio et al.(2020)}]{Barrio-et-al-2020}
		\bibinfo{author}{E.~del Barrio}, \bibinfo{author}{H.~Inouzhe}, \bibinfo{author}{C.~Matrán}, \bibinfo{title}{On approximate validation of models: a Kolmogorov–Smirnov-based approach}, \bibinfo{journal}{TEST}, \bibinfo{volume}{29} (\bibinfo{year}{2020}) \bibinfo{pages}{938--965}.
		
		\bibitem[{Bay and Croix(2017)}]{Bay&Croix2017}
		\bibinfo{author}{X.~Bay}, \bibinfo{author}{J.-C.~Croix}, \bibinfo{title}{Karhunen-Loève decomposition of Gaussian measures on Banach spaces}, \bibinfo{url}{\url{https://arxiv.org/abs/1704.01448}}(\bibinfo{year}{2017}).
		
		\bibitem[{Berlinet and Thomas-Agnan(2011)}]{Berlinet&Thomas-Agnan2011}
		\bibinfo{author}{A.~Berlinet}, \bibinfo{author}{C.~Thomas-Agnan}, \bibinfo{title}{{Reproducing kernel Hilbert spaces in probability and statistics}}, \bibinfo{publisher}{Luxemburg: Springer Sciences and Bussiness Media}, \bibinfo{year}{2011}.
		
		\bibitem[{Biggs et al.(2023)}]{Biggs2023}
		\bibinfo{author}{F.~Biggs}, \bibinfo{author}{A.~Gretton}, \bibinfo{author}{A.~Schrab}, \bibinfo{title}{MMD-FUSE: Learning and Combining Kernels for Two-Sample Testing Without Data Splitting}, \bibinfo{url}{\url{https://arxiv.org/abs/2306.08777}}(\bibinfo{year}{2023}).
		
		\bibitem[{Billingsley(2013)}]{Billingsley2013}
		\bibinfo{author}{P.~Billingsley}, \bibinfo{title}{{Convergence of probability measures}}, \bibinfo{publisher}{John Wiley and Sons}, \bibinfo{year}{2013}.
		
		\bibitem[{Bogachev(1998)}]{Bogachev1998}
		\bibinfo{author}{V.~I.~Bogachev}, \bibinfo{title}{{Gaussian measures}}, \bibinfo{publisher}{Providence: American Mathematical Society}, \bibinfo{year}{1998}.
		
		\bibitem[{Borgwardt et al.(2006)}]{Borgwardt2006}
		\bibinfo{author}{K.~Borgwardt}, \bibinfo{author}{A.~Gretton}, \bibinfo{author}{M.~Rasch}, \bibinfo{author}{B.~Schölkopf}, \bibinfo{author}{A.~J.~Smola}, \bibinfo{title}{A kernel method for the two-sample-problem}, \bibinfo{journal}{Advances in neural information processing systems}, \bibinfo{volume}{19} (\bibinfo{year}{2006}) \bibinfo{pages}{513--520}.
		
		\bibitem[{Borgwardt et al.(2006a)}]{Borgwardt2006a}
		\bibinfo{author}{K.~M.~Borgwardt}, \bibinfo{author}{A.~Gretton}, \bibinfo{author}{H.~P.~Kriegel}, \bibinfo{author}{M.~J.~Rasch}, \bibinfo{author}{B.~Sch\"{o}lkopf}, \bibinfo{author}{A.~J.~Smola}, \bibinfo{title}{Integrating structured biological data by kernel maximum mean discrepancy}, \bibinfo{journal}{Bioinformatics}, \bibinfo{volume}{22} (\bibinfo{year}{2006}) \bibinfo{pages}{e49--e57}.
		
		\bibitem[{Borgwardt et al.(2012)}]{Borgwardt2012}
		\bibinfo{author}{K.~M.~Borgwardt}, \bibinfo{author}{A.~Gretton}, \bibinfo{author}{M.~J.~Rasch}, \bibinfo{author}{B.~Schölkopf}, \bibinfo{author}{A.~Smola}, \bibinfo{title}{A kernel two-sample test}, \bibinfo{journal}{The Journal of Machine Learning Research}, \bibinfo{volume}{13} (\bibinfo{year}{2012}) \bibinfo{pages}{723--773}.
		
		\bibitem[{Brezis(2010)}]{Brezis2010}
		\bibinfo{author}{H.~Brezis}, \bibinfo{title}{{Functional analysis, Sobolev spaces and partial differential equations}}, \bibinfo{publisher}{New York: Springer New York}, \bibinfo{year}{2010}.
		
		\bibitem[{Cárcamo et al.(2020)}]{Carcamo2020}
		\bibinfo{author}{J.~Cárcamo}, \bibinfo{author}{A.~Cuevas}, \bibinfo{author}{L.-A.~Rodríguez}, \bibinfo{title}{Directional differentiability for supremum-type functionals: Statistical applications}, \bibinfo{journal}{Bernoulli}, \bibinfo{volume}{26} (\bibinfo{year}{2020}) \bibinfo{pages}{2143--2175}.
		
		\bibitem[{Cuesta-Albertos et al.(2007)}]{Cuesta2007}
		\bibinfo{author}{J.~A.~Cuesta-Albertos}, \bibinfo{author}{R.~Fraiman}, \bibinfo{author}{T.~Ransford}, \bibinfo{title}{A sharp form of the Cramer–Wold theorem}, \bibinfo{journal}{Journal of Theoretical Probability}, \bibinfo{volume}{20} (\bibinfo{year}{2007}) \bibinfo{pages}{201--209}.
		
		\bibitem[{Cuevas et al.(2004)}]{Cuevas2004}
		\bibinfo{author}{A.~Cuevas}, \bibinfo{author}{M.~Febrero}, \bibinfo{author}{R.~Fraiman}, \bibinfo{title}{An ANOVA test for functional data}, \bibinfo{journal}{Comput. Statist. Data Anal}, \bibinfo{volume}{47} (\bibinfo{year}{2004}) \bibinfo{pages}{111–-122}.
		
		\bibitem[{Debnath and Mikusinski(2005)}]{Debnath&Mikusinski2005}
		\bibinfo{author}{L.~Debnath}, \bibinfo{author}{P.~Mikusinski}, \bibinfo{title}{{Introduction to Hilbert spaces with applications}}, \bibinfo{publisher}{Academic press}, \bibinfo{year}{2005}.
		
		\bibitem[{Dette and Kokot(2022)}]{Dette&Kokot2022}
		\bibinfo{author}{H.~Dette}, \bibinfo{author}{K.~Kokot}, \bibinfo{title}{Detecting relevant differences in the covariance operators of functional time series: A sup-norm approach}, \bibinfo{journal}{Annals of the Institute of Statistical Mathematics}, \bibinfo{volume}{74} (\bibinfo{year}{2022}) \bibinfo{pages}{195--231}.
		
		\bibitem[{Duncan and Wynne(2022)}]{Duncan&Wynne2022}
		\bibinfo{author}{A.~B.~Duncan}, \bibinfo{author}{G.~Wynne}, \bibinfo{title}{A kernel two-sample test for functional data}, \bibinfo{journal}{Journal of Machine Learning Research}, \bibinfo{volume}{23} (\bibinfo{year}{2022}) \bibinfo{pages}{1--51}.
		
		\bibitem[{El-Fallah et al.(2014)}]{El-Fallah2014}
		\bibinfo{author}{O.~El-Fallah}, \bibinfo{author}{K.~Kellay}, \bibinfo{author}{J.~Mashreghi}, \bibinfo{author}{T.~Ransford}, \bibinfo{title}{{A primer on the Dirichlet space}}, \bibinfo{publisher}{Cambridge University Press}, \bibinfo{year}{2014}.
		
		\bibitem[{Eubank and Hsing(2015)}]{Eubank&Hsing2015}
		\bibinfo{author}{R.~Eubank}, \bibinfo{author}{T.~Hsing}, \bibinfo{title}{{Theoretical foundations of functional data analysis, with an introduction to linear operators}}, \bibinfo{publisher}{John Wiley and Sons}, \bibinfo{year}{2015}.
		
		\bibitem[{Fang and Santos(2019)}]{Fang&Santos2019}
		\bibinfo{author}{Z.~Fang}, \bibinfo{author}{A.~Santos}, \bibinfo{title}{Inference on directionally differentiable functions}, \bibinfo{journal}{The Review of Economic Studies}, \bibinfo{volume}{86} (\bibinfo{year}{2019}) \bibinfo{pages}{377--412}.
		
		\bibitem[{Fukumizu et al.(2007)}]{Fukumizu2007}
		\bibinfo{author}{K.~Fukumizu}, \bibinfo{author}{A.~Gretton}, \bibinfo{author}{B.~Schölkopf}, \bibinfo{author}{A.~Smola}, \bibinfo{author}{L.~Song}, \bibinfo{author}{C.~Teo}, \bibinfo{title}{A kernel statistical test of independence}, \bibinfo{journal}{Advances in neural information processing systems}, \bibinfo{volume}{20} (\bibinfo{year}{2007}).
		
		\bibitem[{Fukumizu et al.(2008)}]{Fukumizu2008}
		\bibinfo{author}{K.~Fukumizu}, \bibinfo{author}{A.~Gretton}, \bibinfo{author}{B.~Sch\"{o}lkopf}, \bibinfo{author}{X.~Sun}, \bibinfo{title}{Kernel measures of conditional dependence}, \bibinfo{journal}{Advances in neural information processing systems}, \bibinfo{volume}{20} (\bibinfo{year}{2008}) \bibinfo{pages}{489--496}.
		
		\bibitem[{Fukumizu et al.(2009)}]{Fukumizu2009}
		\bibinfo{author}{K.~Fukumizu}, \bibinfo{author}{A.~Gretton}, \bibinfo{author}{G.~R.~Lanckriet}, \bibinfo{author}{B.~Sch\"{o}lkopf}, \bibinfo{author}{B.~K.~Sriperumbudur}, \bibinfo{title}{Kernel choice and classifiability for RKHS embeddings of probability distributions}, \bibinfo{journal}{Advances in neural information processing systems}, \bibinfo{volume}{22} (\bibinfo{year}{2009}) \bibinfo{pages}{}.
		
		\bibitem[{Fukumizu et al.(2010)}]{Fukumizu2010}
		\bibinfo{author}{K.~Fukumizu}, \bibinfo{author}{A.~Gretton}, \bibinfo{author}{G.~R.~Lanckriet}, \bibinfo{author}{B.~Sch\"{o}lkopf}, \bibinfo{author}{B.~K.~Sriperumbudur}, \bibinfo{title}{Hilbert space embeddings and metrics on probability measures}, \bibinfo{journal}{The Journal of Machine Learning Research}, \bibinfo{volume}{11} (\bibinfo{year}{2010}) \bibinfo{pages}{1517--1561}.
		
		\bibitem[{Fukumizu et al.(2011)}]{Fukumizu2011}
		\bibinfo{author}{K.~Fukumizu}, \bibinfo{author}{G.~R.~Lanckriet}, \bibinfo{author}{B.~K.~Sriperumbudur}, \bibinfo{title}{Universality, characteristic kernels and RKHS embedding of measures}, \bibinfo{journal}{Journal of Machine Learning Research}, \bibinfo{volume}{12} (\bibinfo{year}{2011}) \bibinfo{pages}{2389--2410}.
		
		\bibitem[{Fukumizu et al.(2013)}]{Fukumizu2013}
		\bibinfo{author}{K.~Fukumizu}, \bibinfo{author}{A.~Gretton}, \bibinfo{author}{D.~Sejdinovic}, \bibinfo{author}{B.~Sriperumbudur}, \bibinfo{title}{Equivalence of distance-based and RKHS-based statistics in hypothesis testing}, \bibinfo{journal}{The annals of statistics}, \bibinfo{volume}{} (\bibinfo{year}{2013}) \bibinfo{pages}{2263-2291}.
		
		\bibitem[{Ghosh et al.(2016)}]{Ghosh2016}
		\bibinfo{author}{S.~Ghosh}, \bibinfo{author}{G.~M.~Pomann}, \bibinfo{author}{A.~M.~Staicu}, \bibinfo{title}{A two‐sample distribution‐free test for functional data with application to a diffusion tensor imaging study of multiple sclerosis}, \bibinfo{journal}{Journal of the Royal Statistical Society, Ser. C}, \bibinfo{volume}{65} (\bibinfo{year}{2016}) \bibinfo{pages}{395--414}.
		
		\bibitem[{Giné and Nickl(2008)}]{Gine&Nickl2008}
		\bibinfo{author}{E.~Giné}, \bibinfo{author}{R.~Nickl}, \bibinfo{title}{Uniform central limit theorems for kernel density estimators}, \bibinfo{journal}{Probability Theory and Related Fields}, \bibinfo{volume}{141} (\bibinfo{year}{2008}) \bibinfo{pages}{333--387}.
		
		\bibitem[{Giné and Nickl(2016)}]{Gine&Nickl2016}
		\bibinfo{author}{E.~Giné}, \bibinfo{author}{R.~Nickl}, \bibinfo{title}{{Mathematical foundations of infinite-dimensional statistical models}}, \bibinfo{publisher}{Cambridge Series in Statistical and Probabilistic Mathematics. Cambridge University Press}, \bibinfo{year}{2016}.
		
		\bibitem[{Gretton et al.(2020)}]{Gretton2020}
		\bibinfo{author}{A.~Gretton}, \bibinfo{author}{F.~Liu}, \bibinfo{author}{J.~Lu}, \bibinfo{author}{D.~J.~Sutherland}, \bibinfo{author}{B.~Laurent}, \bibinfo{author}{W.~Zhang}, \bibinfo{title}{Learning deep kernels for non-parametric two-sample tests}, \bibinfo{journal}{International conference on machine learning}, \bibinfo{volume}{} (\bibinfo{year}{2020}) \bibinfo{pages}{6316--6326}.
		
		\bibitem[{Guo et al.(2022)}]{Guo2022}
		\bibinfo{author}{J.~Guo}, \bibinfo{author}{J.~T.~Zhang}, \bibinfo{author}{B.~Zhou}, \bibinfo{title}{Testing equality of several distributions in separable metric spaces: A maximum mean discrepancy based approach}, \bibinfo{journal}{Journal of Econometrics}, \bibinfo{volume}{} (\bibinfo{year}{2022}) \bibinfo{pages}{}.
		
		\bibitem[{Hall and Van Keilegom(2007)}]{Hall2007}
		\bibinfo{author}{P.~Hall}, \bibinfo{author}{I.~Van Keilegom}, \bibinfo{title}{Two-sample tests in functional data analysis starting from discrete data}, \bibinfo{journal}{Statistica Sinica}, \bibinfo{volume}{17} (\bibinfo{year}{2007}) \bibinfo{pages}{151--1531}.
		
		\bibitem[{Janson(1997)}]{Janson1997}
		\bibinfo{author}{S.~Janson}, \bibinfo{title}{{Gaussian Hilbert Spaces}}, \bibinfo{publisher}{Cambridge University Press}, \bibinfo{year}{1997}.
		
		\bibitem[{Ledoux and Talagrand(2013)}]{Ledoux&Talagrand2013}
		\bibinfo{author}{M.~Ledoux}, \bibinfo{author}{M.~Talagrand}, \bibinfo{title}{{Probability in Banach Spaces: Isoperimetry and processes}}, \bibinfo{publisher}{Springer Science and Business Media}, \bibinfo{year}{2013}.
		
	\bibitem[{Lehmann and Romano(2022)}]{Lehmann&Romano2022}
		\bibinfo{author}{E.~L.~Lehmann}, \bibinfo{author}{J.~P.~Romano}, \bibinfo{title}{{Testing Statistical Hypotheses}}, \bibinfo{publisher}{Springer Texts in Statistics}, \bibinfo{year}{2022}.



		\bibitem[{Marcus(1985)}]{Marcus1985}
		\bibinfo{author}{D.~J.~Marcus}, \bibinfo{title}{Relationships between Donsker classes and Sobolev spaces}, \bibinfo{journal}{Zeitschrift für Wahrscheinlichkeitstheorie und Verwandte Gebiete}, \bibinfo{volume}{69} (\bibinfo{year}{1985}) \bibinfo{pages}{323--330}.
		
		\bibitem[{Müller(1997)}]{Muller1997}
		\bibinfo{author}{A.~Müller}, \bibinfo{title}{Integral probability metrics and their generating classes of functions}, \bibinfo{journal}{Advances in Applied Probability}, \bibinfo{volume}{29} (\bibinfo{year}{1997}) \bibinfo{pages}{429--443}.
		
		\bibitem[{Paulsen and Raghupathi(2016)}]{Paulsen2016}
		\bibinfo{author}{V.~I.~Paulsen}, \bibinfo{author}{M.~Raghupathi}, \bibinfo{title}{{An introduction to the theory of reproducing kernel Hilbert spaces}}, \bibinfo{publisher}{Cambridge University Press}, \bibinfo{year}{2016}.
		
		\bibitem[{Rizzo and Sz\'{e}kely(2017)}]{Rizzo&Szekely2017}
		\bibinfo{author}{M.~L.~Rizzo}, \bibinfo{author}{G.~J.~Sz\'{e}kely}, \bibinfo{title}{The energy of data}, \bibinfo{journal}{Annual Review of Statistics and Its Application}, \bibinfo{volume}{4} (\bibinfo{year}{2017}) \bibinfo{pages}{447--479}.
		
		\bibitem[{Rizzo and Szekely(2022)}]{EnergyPackage2022}
		\bibinfo{author}{M.~Rizzo}, \bibinfo{author}{G.~Szekely}, \bibinfo{title}{E-statistics: Multivariate inference via the energy of data. R package version 1.7--10}, \bibinfo{url}{\url{https://CRAN.R-project.org/package=energy}}, (\bibinfo{year}{2022}).
		
		\bibitem[{Scholk\"{o}pf and Smola(2018)}]{Scholkopf&Smola2018}
		\bibinfo{author}{B.~Sch\"{o}lkopf}, \bibinfo{author}{A.~J.~Smola}, \bibinfo{title}{{Learning with kernels: Support vector machines, regularization, optimization, and beyond}}, \bibinfo{publisher}{MIT press}, \bibinfo{year}{2018}.
		
		\bibitem[{Shapiro(1991)}]{Shapiro1991}
		\bibinfo{author}{A.~Shapiro}, \bibinfo{title}{Asymptotic analysis of stochastic programs}, \bibinfo{journal}{Annals of Operations Research}, \bibinfo{volume}{30} (\bibinfo{year}{1991}) \bibinfo{pages}{169--186}.
		
		\bibitem[{Smaga and Zhang(2022)}]{Smaga&Zhang2022}
		\bibinfo{author}{L.~Smaga}, \bibinfo{author}{J.~T.~Zhang}, \bibinfo{title}{Two-sample test for equal distributions in separate metric space: New maximum mean discrepancy based approaches}, \bibinfo{journal}{Electronic Journal of Statistics}, \bibinfo{volume}{16} (\bibinfo{year}{2022}) \bibinfo{pages}{4090--4132}.
		
		\bibitem[{Sriperumbudur(2016)}]{Sriperumbudur2016}
		\bibinfo{author}{B.~Sriperumbudur}, \bibinfo{title}{On the optimal estimation of probability measures in weak and strong topologies}, \bibinfo{journal}{Bernoulli}, \bibinfo{volume}{22} (\bibinfo{year}{2016}) \bibinfo{pages}{1839--1893}.
		
		\bibitem[{Tukey(1959)}]{Tukey1959}
		\bibinfo{author}{J.~W.~Tukey}, \bibinfo{title}{A quick compact two sample test to Duckworth's specifications}, \bibinfo{journal}{Technometrics}, \bibinfo{volume}{1} (\bibinfo{year}{1959}) \bibinfo{pages}{31--48}.
		
		\bibitem[{van der Vaart and Wellner(2023)}]{van_der_Vaart&Wellner2023}
		\bibinfo{author}{A.~W.~van der Vaart}, \bibinfo{author}{J.~A.~Wellner}, \bibinfo{title}{{Weak convergence and empirical processes: With applications to statistics}}, \bibinfo{publisher}{Springer}, \bibinfo{year}{2023}.
		
		\bibitem[{Vapnik(2013)}]{Vapnik2013}
		\bibinfo{author}{V.~Vapnik}, \bibinfo{title}{{The nature of statistical learning theory}}, \bibinfo{publisher}{Springer science and business media}, \bibinfo{year}{2013}.
		
		\bibitem[{Wendland(2004)}]{Wendland2004}
		\bibinfo{author}{H.~Wendland}, \bibinfo{title}{{Scattered data approximation}}, \bibinfo{publisher}{Cambridge university press}, \bibinfo{year}{2004}.
		
		\bibitem[{Zhang and Zhao(2013)}]{ZhangZhao2013}
		\bibinfo{author}{H.~Zhang}, \bibinfo{author}{L.~Zhao}, \bibinfo{title}{On the inclusion relation of reproducing kernel Hilbert spaces}, \bibinfo{journal}{Analysis and Applications}, \bibinfo{volume}{11} (\bibinfo{year}{2013}) \bibinfo{pages}{}.
		
		\bibitem[{Zhang(2014)}]{Zhang2014}
		\bibinfo{author}{J.~T.~Zhang}, \bibinfo{title}{{Analysis of variance for functional data}}, \bibinfo{publisher}{CRC Press}, \bibinfo{year}{2014}.
	\end{thebibliography}
		%
		%

\medskip

\end{document}